\documentclass[11pt, a4paper, openany]{amsart}
\usepackage[applemac]{inputenc}  
\usepackage{t1enc}    
\usepackage[all,cmtip,2cell]{xy}
\UseTwocells
\usepackage{ucs} 
\usepackage[applemac]{inputenc}
\usepackage{textcomp}  
\usepackage{t1enc}    
\usepackage{subfiles}
\usepackage{article}
\usepackage{endnotes}
\usepackage{wasysym}
\usepackage{ifsym}
\usepackage[mathcal]{eucal}
\usepackage[  hidelinks,
  bookmarks=true,
  bookmarksnumbered=true,
  pdfusetitle]{hyperref} 
\usepackage{mathrsfs}
\usepackage{yfonts}
\usepackage{stmaryrd}
\usepackage{xcolor}
\usepackage{lmodern}
\usepackage{slantsc}

\DeclareMathOperator{\Epi}{Epi}
\DeclareMathOperator{\dom}{dom}

\DeclareMathOperator{\op}{op}
\newcommand{\Si}{\mathfrak{S}}
\newcommand{\true}{\mathrm{true}}
\newcommand{\false}{\mathrm{false}}
\newcommand{\U}{\mathcal{U}}

\hypersetup{
  pdftitle    = {Bridging Badiou's Gaps: Topos Theory and the Bodies of Mathematical Expression},
  pdfauthor   = {Antti Veilahti},
  pdfsubject  = {Topos theory, sheaf theory, philosophy of mathematics},
  pdfkeywords = {topos theory, Badiou, forcing, sheaves, Heyting algebra},
  pdfcreator  = {LaTeX with hyperref}
}

\usepackage{fancyhdr}
\usepackage{ifthen}
\author{Antti Veilahti \vspace{.5cm}}

\title{\Large{Bridging Badiou’s Gaps: \\ $\ $ \\ Topos Theory and the Bodies of Mathematical Expression }}

\begin{document}

\begin{abstract}

This article develops a second stage of an argument begun in 2013, revisiting Alain Badiou's use of topos theory in \textit{Being and Event} and \textit{Logics of Worlds}. The previous article 
diagnosed Badiou's reduction of topos theory to its atomistic, localic form as the scientific scapegoat in his attack on `democratic materialism' \href{https://arxiv.org/pdf/1301.1203}{(arXiv:1301.1203)}. Now I ask what follows when this logicist illusion is dispelled and the categorical alternative is taken seriously, particularly with respect to forcing, semantic expression, and subjectification that constitute the key elements of Badiou's earlier work.

 I reconstruct Badiou's own arguments in a corrected form, while transferring them to a categorical framework.  The aim is not to introduce new mathematical results, but to identify what is conceptually lost through their reduction to set-theoretic formalism. In particular, I describe Cohen's forcing argument in the context of Grothendieck topoi and Kripke--Joyal semantics, thereby explicating the geometric structures implicit in Badiou's thinking.

This reconstruction further brings into view a different way of understanding Badiou's concepts of subject, body, and change. Rather than treating a topos as merely a formal container for mathematical objects, I propose reading it as a body of mathematical \textit{expression}, aligning with Johnstone's characterization of a topos as a `universe of discourse'. At the same time, the geometrization of logic allows Badiou's own four-fold typology of change---when applied to topos theory itself as a mathematical `event'---to be anchored to various historical definitions of a topos. 

In this way, the article bridges Badiou's two major works while clarifying how his own limited expression of mathematical formalism embodies his `dialectical', ultimately singular concept of the event, and the way it clouds his philosophical understanding of locality and temporality.

\end{abstract}
\vspace{-1cm}
\maketitle	
\thispagestyle{empty} 
\pagenumbering{roman}

\tableofcontents



\pagestyle{plain}
\pagenumbering{arabic}

\makeatletter
\def\@makefnmark{\hbox{\@textsuperscript{\normalfont\@thefnmark}}}
\makeatother

\section*{Introduction}

Contemporary physics---the science of matter---tackles highly refined and technically complex theories, leading into its `Church' becoming increasingly divided, to use a Hegelian metaphor. However, what is divided is less matter itself than our ways of approaching this division. As Nietzsche polemically declared: `we must [\dots] declare war---a ruthless fight' not primarily against the atomism of matter, but `[f]irst of all [\dots] that other and more disastrous atomism [\dots] the \textit{atomism of the soul}'\footnote{See Nietzsche, Friedrich (2002), \textit{Beyond Good and Evil}, translated by Judith Norman, Cambridge University Press, p. 14.}.

It is this intellectual diaspora, moving away from simplistic, atomistic representations, that appears to bring a `dialectical' momentum into science---i.e. an `infinite truth-process' that Alain Badiou\footnote{See Badiou, Alain (2006), \textit{Being and Event}, translated by Oliver Feltham, Continuum.} situated alongside `art', `politics', and `love'---like my first article explored\footnote{The first article of this two-stage intervention is the substantially revised 2026 version of Veilahti, Antti (2013), `Alain Badiou's Mistake---Two Postulates of Dialectical Materialism', \textit{Preprint}, arXiv. \url{https://arxiv.org/pdf/1301.1203}.}. In particular, the changing process of scientific expression appears to contradict the Peircean understanding of science as a limit-process that would eventually neutralize an otherwise infinite Truth. To Badiou, such neutrality amounts to a deception; rather, truth is always at \textit{risk}. Crucially, this risk does not concern only the empirical sciences: it is a risk inherent to the very science of being itself---that is, mathematical ontology. 

Physics, too, increasingly puts at risk the conventional modes of scientific expression, as I also previously examined\footnote{See Veilahti, Antti (2017), `Higher Theory and the Three Problems of Physics', \textit{Preprint}, arXiv. \url{https://arxiv.org/pdf/1712.09454}.}. Moreover, category theory, and topos theory as one of its most philosophically compelling applications, have provided a mathematical framework in which the limits of scientific expression themselves become objects of inquiry. 
	
Badiou, by contrast, believes this risk itself to lie outside the domain of scientific scrutiny. \textit{Being and Event}---which rendered Badiou one of the most prominent philosophers in France upon its publication in 1988---could thus only represent this risk indirectly, by \textit{locating} it within an evental, `paradoxical' multiple $e \in e$: something formally decoupled from scientific expression. In fact, if we follow his doctrine that mathematics is ontology, then the Hilbertian program would carry the mantle of `meta-ontology' as the proper form of scientific subjectivity. The only way for Badiou to objectify the event formally was therefore to assign it negatively through the supposedly `paradoxical' language of set theory, thus effectively presenting the event through double negation.

This formal language then became the only, even if internally contingent, regulatory framework through which to fully appreciate the event as it breaks the axiom of foundation\footnote{\label{axiomoffoundation}The axiom of foundation states that for any non-empty set $x$, there exists an element $y \in x$ such that $x \cap y = \emp$.}. When it comes to the philosophical momentum of his own endeavor, it was thus this \textit{language} itself, anomalous as it may seem, which became `another event'---a distinct form of `fidelity' separated from the event itself by a gap that Badiou himself labeled as its `\textit{durée}'.
And since an event is always a risk, it was not surprising that subjectivity inherent to mathematical inquiry itself becomes visible through a procedure established by Paul J. Cohen;\footnote{See Cohen, Paul J. (1963), `The Independence of the Continuum Hypothesis', \textit{Proc. Natl. Acad. Sci. USA}, 50(6), pp. 1143--1148; and Cohen, Paul J. (1964), `The Independence of the Continuum Hypothesis II', \textit{Proc. Natl. Acad. Sci. USA}, 51(1), pp. 105--110.} Cohen's argument produces a so-called `generic filter' ($\female$) that captures the self-referentiality inherent to logical inquiry itself, thus establishing an analogue between the ontological event $e$ and the formal regulation of being expressed in mathematical language, treated thus as `another event'.

In effect, `\textit{durée}' is not just an empty phrase separating the event $e$ and the sign of intervention $\square$; instead, it is Badiou's own \textit{passage} from the fidelity of the event---the first half of \BE---toward the tracing of the `generic' procedure that elucidates mathematical inquiry itself as a temporal process. Therefore, while framing the duration of the event as a `gap' might be a formal outcome of the specific intervention itself, this does not lessen the significance of this formalism as a core structure in Badiou's philosophy, where it marks a continuity between the two aspects of the same anomalous encounter: the event itself and its local and subjective appearance by the means of mathematical expression. 

Rather than depicting the event as such, \BE\ thus mirrors an antagonism inherent to the event itself by breeding a philosophical gap between the `meta-ontological' status of set theory---as a supposedly causal rupture---and the true consequences of this illogical rift for Badiou's formalist, `dogmatic fidelity'; consequences made apparent and yet structurally underexplored by Cohen's model-theoretic observations. Thus, \BE\ failed to capture the precise meaning of this `gap' in a way consequential to the formalist core of Badiou's own philosophy, that is, as it overwhelms his own, quasi-formalist expression. Insofar as Cohen's procedure is concerned, and to endorse Badiou's own Dieudonné-inspired maxim of the need to master the contemporary mathematical corpus\footnote{See Badiou, Alain (2006), \textit{Briefings on Existence: A Short Treatise on Transitory Ontology}, translated by Norman Madarasz, State University of New York Press, p. 12.}, my paper claims that the answer must be structurally located within \textit{topos theory}. 

Indeed, soon after Cohen had established the independence of the continuum hypothesis, William Lawvere and Myles Tierney recognized a structural analogue between this procedure and the categorically oriented methods Alexander Grothendieck had applied to algebraic geometry. This led to the development of topos theory which, distancing itself from the direct, `analytic' representation of logic, was able to produce a set-theoretic (Mitchell--Bénabou) language relationally, without explicitly forcing out that geometric torsion which set theory violently denies by its binary, flat representation of semantics. In other words, semantics itself became a place for geometrically non-trivial manifestations---that is, an admission of the \textit{real} or \textit{material} possibility of subjective, truth-related torsion. 

Failing to bridge with this mathematical passage, the `gap'---or temporal aspect of the event---remains ultimately neutralized in Badiou's own philosophical endeavor: this fact is made evident by the absence of a true disconnect between his two works. In other words, while recognizing that `problems of connection and continuity'---which he leaves `for others to solve'---remain `between `generic procedure’ and `intra-worldly consequences of the existence of an inexistent’\footnote{Badiou, \LW, 2009, p. 39.}, this technical gap conceals the more pervasive, meta-ontological gap between the two virtually distinct modes of formalism the works engage. It is the appreciation of this second, ontological antagonism that is the secret to conceiving mathematical expression itself as an event---an intervention that, indeed, is ultimately both temporal and historical. 

Thus his formally subordinate mathematical posture marks rather a stasis that conceals the true `gap' concerning \textit{mathematical expression} itself. Symptomatically, Badiou never explicates the distance between Cohen's and Russell's insights, or, as he was later made aware by `perspicacious readers'---namely `Desanti, Deleuze, Nancy and Lyotard' as he himself mentions---the fact that the framing of the ontological definition of `what happens' occurs both from below, by `positing the existence [\dots] of an event-site', and from above, by `demanding that every event receive a name'\footnote{Badiou, \LW, 2009. p. 361.}. 
Despite seeking to address these concerns in \LW, Badiou\footnote{Ibid., p. 357.} himself confesses how both works were haunted by certain negativity, as `we cannot find the means to identify change either in the order of mathematics [\dots] or in the logic, the thinking of being-there or appearing'. Even if he had `not taken care to guarantee at every point a continuity between the two projects', both works rely on the same premises of a `pure theory of a multiple', grounded on his belief that `some ontological reminders to suffice as the guarantee of [their] homogeneity'\footnote{Ibid., p. 39.}.

As a consequence, much of \LW\ remains reactionary in its attempt to establish a `real synthesis' between the `ontological' and `phenomenological' realms---an attempt that my previous article\footnote{See Veilahti, `Alain Badiou's Mistake', 2013.} demonstrated as being atomistic, `caught up in the transcendental indexing of which it is the ontological support'\footnote{Badiou, Alain, \textit{Logics of Worlds}, translated by Alberto Toscano, Continuum, 2009, p. 360.}, instead of marking a true breakage to understand the individuation of mathematical expression itself. His deliberate effort to reduce mathematics to the mold of set-theory and localic topos theory narrows it to the extent that precludes any more dynamic, process philosophical interpretation, let alone a novel vein of political thinking. Therefore, what connects Badiou's two works is actually an intellectual disconnect that was immanent already in \BE's original struggle to envision a subject, or intervention, in a way capable of expressing subjective torsion at the most material, formal level---not merely as a paradox of self-reference. 

Nonetheless, the seeds of a genuinely \textit{topological} and geometric understanding of the event were already forming in \BE, and their adequate structural articulation could have made \LW\ a genuine `phenomenological' contribution to event-philosophy, surpassing its logic-bound composition. Instead of internalizing the struggle by focusing on the `intra-worldly consequences of the existence of an inexistent', like Badiou formally does, this \textit{philosophical expansion} requires establishing how precisely Cohen's procedure relates to the structural ambiguity and torsion inherent to \textit{localization} as a concept---i.e., the very meaning of the distance between what happens `from below' and how it is nominated `from above'. 
Presenting this philosophical controversy merely as a formal `gap' specifically masks its momentum and impedes our grasp of \textit{how the concept of locality can itself be contextualized}, thus treating the mathematical procedure or expression itself as an event to be located. This, in other words, would illustrate the true meaning of the `deductive fidelity' of science, or mathematical expression, counting as another event.

Badiou, in fact, seems to have recognized the value of this conceptual shift intuitively, by directly engaging with the concept of a topos in \LW, though retreating into his established formalism as he replaces the general notion of a topos solely by the localic category $\sets^{T^{\op}}$ for an explicit transcendental indexing lattice $T$, as my previous paper established\footnote{See Veilahti, `Alain Badiou's Mistake', 2013.}. This reduction implies a specific restriction on the \textit{form} of the subjectification of truth---a `point of truth'---like this paper explores. Namely, Badiou can formally conceive only a \textit{split} choice (either/or), instead of grasping the more, in Kierkegaardian terms, `radical' choice invoked by category theory. Badiou thus continues to overlook the `diagrammatic' paradigm behind the categorical mode of expression, in spite of its historical relevance in the works of figures such as Alexander Grothendieck and Pierre Deligne\footnote{For a comprehensive philosophical account of the categorical shift, see Krömer, Ralf (2007), \textit{Tool and Object: A History and Philosophy of Category Theory}, Birkhäuser; see also Marquis, Jean-Pierre (2009), \textit{From a Geometrical Point of View: A Study of the History and Philosophy of Category Theory}, Springer.}.

This logical asylum persists in the third iteration of the \textit{Being and Event} trilogy. In \textit{Immanence of Truths}\footnote{Badiou, Alain (2022), \textit{The Immanence of Truths: Being and Event III}, translated by A. Toscano and K. Reinhard, Bloomsbury Academic.}, Badiou's recourse to large cardinals and transfinite analytics under the arithmetic boundary of $\kappa$-compactness again attempts to save the `absolute' from the prolific multiplicity of mathematical expressions, treating truth as a static creation `outside of knowledge' and as a mere dialectical negative of his local, hierarchical concept of time\footnote{For an alternative framework that embraces this multiplicity not as an ontological threat but as an interactive, relational starting point, also see Zalamea, Fernando (2014), \textit{Synthetic Philosophy of Contemporary Mathematics}, translated by Z. L. Fraser, Urbanomic.}.

My forthcoming philosophical autoethnography will examine this maneuver in detail, arguing that the compactness Badiou elevates into an ontological principle finds a revealing categorical analogue in the accessibility conditions through which Lurie---as a mathematician, or an `unconscious ontologist' as Badiou calls them---enforces the ontological tractability of higher topoi. Literally speaking, of course, the occurrence of `$\kappa$-compactness' is accidental as Lurie's and Badiou's notions are mathematically distinct, reflecting different approaches to actualization---arithmetic and geometric, respectively.  Yet their juxtaposition exposes a common philosophical problem: how an apparently indefinite multiplicity is conditioned (and forced) to become formally accessible. 

This makes Lurie's project resonate paradigmatically with Badiou's `dialectical' project. The relevance of the latter does not reside within the virtuality of higher homotopies \textit{per se}; rather, it emerges through the problem of its domestication through specific, ontologically controlled localization functors. Even if `the same mathematics' were to still prevail, it would only do so by continuing to impose similar laws of subjectivity depicted in \BE, again equipping the subject with essentially only two options: either to abandon formal ontology entirely, or to accept it as the only scientifically legitimate basis of access. 

\textit{Being and Event III} might have had something to say about this process, especially given how I personally informed Badiou about my project years ago. Unfortunately, it recedes from the non-local, interactional geometry that both \BE\ and \LW\ flirted with, instead of contributing to a true bridge between Badiou's work and geometry. \textit{Immanence of Truths}, which in its own words defines geometry as a withering `system of relationships that an object maintains with its context', is still preoccupied by the question of totalities that `geometrize[] all thought', never moving beyond this dialectical process immanent to mathematical ontology. It seeks to `know' those higher bodies from the outset, instead of posing the speculative, but at the same time empiricist question of what things of the world, which supposedly embody such higher structures\footnote{See Veilahti, `Higher Theory and the Three Problems of Physics', 2017.}, `know' or experience within themselves.

Nonetheless, in the present article I demonstrate that the intuitive shift between Badiou's first two masterpieces, when rephrased in category-\bk theoretical terms, retains the potential to overcome the partialism inherent to Badiou's own formalism---a new `style of life’ that Badiou himself, as a philosophical persona, does not yet recognize. Therefore, where Badiou discusses the `materialism' of the object, whether atomistic or truly dialectical, topos theory proper allows an understanding of the \textit{materialism of the subject}---that is, how truth itself emerges as a local and localizing effect.

This consists of a relational rereading not just of what mathematics expresses, but also of \textit{how} it expresses and deduces the `points' articulated in formal language, shifting the focus away from the contents of objects to their different modes of expression. As such, topos theory turns into a magnificent `boundary object'\footnote{See Bowker, Geoffrey C. and Susan Leigh Star (1999), \textit{Sorting Things Out: Classification and its Consequences}, MIT Press.} or bridge---not just between philosophy and science, or objectivity and subjectivity, but contrasting entire schools of philosophical reasoning starting from Parmenides and Heraclitus, one occupied by the `dialectical' problems of being and one inclined towards concepts of difference and becoming. 

This elaboration is not merely a technical nuance in Badiou scholarship; it bears broader implications to the philosophy of mathematics. The central thesis is my framing of mathematical immanence, or what Badiou labels as `Truths', as a historical process. This is not a history deployed in anthropocentric coordinates. Rather, it is one in which mathematical abstraction and innovation themselves unfold temporally. The concept of a topos gives this process a mathematical expression: a site that is speculative yet empirical in the sense that it allows experimenting with mathematical constructs\footnote{The experimental aspect of topos theory is essentially twofold:  its empirical value emerges both through its ability to express mathematical innovation and change, as well as through its implications to physics, as I have reviewed in Veilahti, `Higher Theory and the Three Problems of Physics', 2017.}. 

In the present article, the `singular' 1-topos remains the principal formal setting. Read through my topos-theoretic reconstruction of Badiou's fourfold concept of change,  it resonates with Hegel's dialectical conception of history. By contrast, my forthcoming work turns to higher topoi while endorsing an asingular and interactive approach to mathematical temporality, thereby evoking a concept of history that is multiplicative and differential rather than dialectical. 

This review of Badiou's work, which combines a formal, mathematical elaboration with a profound analysis of its implications for Badiou's philosophical project, is organized as follows. Drawing on the schematic approach already deployed in the first article, Sections 1--2 examine the question of topology from two alternative points of view: either by identifying the properly relational conditions that ensure the elicitability of ontological hypotheses, or by proceeding through Badiou's `inside-out' path where set-theoretic boundedness is taken for granted, attempting to reconstruct `diagrammatic' patterns though with no regulatory value. There, where the first article focused on the question of materiality, Section 2 focuses instead on points as a dualistic couple, illustrating the conditions on which also a categorical topos, when subjected to an ontological gaze emanating from a `Platonic' one, is forced to answer where, and through what point, it makes a difference. 

Sections 3--5 then transfer this problematic to the context of Cohen's procedure, again overturning the reading of the procedure in \BE, while demonstrating the structural, categorical precursors that effectuate how Cohen's procedure ultimately emerges as a topological moment. This culminates in the actual independence proof, which at one and the same time both corrects shortcomings inherent to Badiou's work while demonstrating conceptually how the relational and geometric approach to points and subjectivity transform Cohen's results into a paradigm, not an anomaly. 

Sections 6--7 complete this mathematical portrait by fully reverting what is left of Badiou's formalist program: through sheaf theory, it can be illustrated how logical semantics and ultimately Badiou's concept of mathematical intervention can be conceived as products \textit{internal} to a topos---that is, as entities born from the body of mathematical expression, rather than holding it together. Beyond sheaf theory, which is still `dialectically materialist' in Badiou's sense of atomistic ontology, a genuinely categorical approach to Kripke--Joyal semantics demonstrates how expressions of truth and subjectivity ultimately require a contextual understanding, and where Badiou's philosophy can thus be understood as philosophy inside a topos, rather than about it. 

Sections 8--9, once more, transfer the implications drawn from these mathematical observations to the philosophical domain. These sections are written so that they can be read independently of the mathematically oriented Sections 1--7. The particular focus is on the relationship between the event and mathematics, as Badiou at times conceives them as absolutely distinct---
thus infecting the politics of the event with entitlement and grandiosity---while at other times calls mathematical intervention itself `another event': a picture particularly completed by the results from Section 7 that permit a topological understanding of a multiplicity of mathematical interventions and temporalities. The last section further connects this analysis to questions of change and a body---concepts pivotal to the more phenomenologically oriented work of \LW---ultimately arguing that elementary topos theory is a plausible framework for a formal reinterpretation of Badiou's `dialectical' philosophy but that in itself it falls short of more post-structuralist orientations where topos theory itself is not to be conceived as a fixed, singular outcome but rather a process.

The main works through which I re-read Badiou's still dialectical philosophical trajectory are the classical foundations by Johnstone\footnote{See Johnstone, Peter T. (1977), \textit{Topos Theory}, Academic Press; see also Johnstone, Peter T. (2002), \textit{Sketches of an Elephant: A Topos Theory Compendium}, Volume 1, Clarendon Press.} and Mac Lane and Moerdijk\footnote{See Mac Lane, Saunders and Ieke Moerdijk (1992), \textit{Sheaves in Geometry and Logic: A First Introduction to Topos Theory}, Springer-Verlag.}, intersected with contemporary reflections on the conceptual shift toward higher topologies\footnote{See, e.g., Anel, Mathieu and Gabriel Catren (eds.) (2022), \textit{New Spaces in Mathematics: Formal and Conceptual Reflections}, Cambridge University Press.}. 
The mathematical results employed here are standard; what is at issue is their reconstruction and categorical reformulation within arguments that until now have remained philosophically constrained by Badiou's own, fully set-theoretic presentation.

\section{From Localic Frames to the Categorical Expression of Torsion }

In a previous article\footnote{Veilahti, `Alain Badiou's Mistake', 2013.}, I reframed Badiou's\footnote{Badiou, \LW, 2009.} recent approach to `calculated phenomenology' in terms of categorical topos theory. It turns out that Badiou grounds his philosophical argument on a restricted understanding now recognized as the regime of `locales'---the so-called \textit{local} Grothendieck topoi, or $0$-topoi as they are called in higher topos theory\footnote{Lurie, Jacob (2009), \textit{Higher Topos Theory}. Princeton: Princeton University Press.}. These are mathematical structures graded on an \textit{external complete Heyting algebra}; that is, a lattice $T$ with minimal and maximal elements, a distributive symmetric conjunction $\wedge: T \times T \to T$, and a completion operator $\Sigma: \Power T \to T$ which marks the least upper bound $\Sigma S$ for a subset $S \subset T$. 

This architecture gives rise to a logic `internal' to this \textit{external} algebra, where $\wedge$ denotes conjunction, implication is defined as 
$$p \Rightarrow q = \Sigma \{t \mid p \wedge t \leq q\},$$
and negation is given by 
$$ \neg p = \Sigma \{q \mid p \wedge q = \mu\}\footnote{Badiou utilizes the Greek letter $\mu$ to designate the minimal degree (\textit{le degré minimal}) of a transcendental framework, denoting that which `inexists' within a specific world (see Badiou, \LW, 2009, pp. 119--121). Formally, in the context of an external complete Heyting algebra, $\mu$ operates as the bottom element ($\bot$), functioning as the structural zero of the localic lattice. I retain Badiou's notation to preserve the ontological status of the inexistent, which is structurally localized rather than mathematically void.}.$$
One could say that $T$ measures the `degrees' of truth, while the maximal and minimal elements refer to absolute truth and falsehood. In particular, the logic internal to this external algebra satisfies the axioms of intuitionist logic, where double negation $\neg\neg$ could be non-trivial. By contrast, in a general elementary topos, logic emerges as a relational expression of the so-called \textit{internal} complete Heyting algebra $\Omega$, taking the form of the Mitchell--Bénabou language and Kripke--Joyal semantics, as I will later illustrate.

It was Badiou's explicit aim to move toward category theory to account for certain questions regarding Kant's suggestion of the \textit{a priori} forms of sensibility and transphenomenal real\footnote{Badiou, \LW, 2009, p. 104.}. This transition gives rise to the designation of an `object' approached relationally rather than through absolute, incorporeal `contents', as if to protect it from Kant's `legalism---always asking \textit{Quid juris?} or ``Haven’t you crossed the limit?'' ` Within this framework, an object is no longer defined by elements it inherently consists of or incorporates, but strictly by how it relates to other objects through the so-called arrows or morphisms.
 
Yet, despite his Kantian remarks, Badiou never espoused this transition to categorical mathematics formally, even in his latest iteration of the \textit{Being and Event} -series. In particular, he attempts to formulate the underlying mathematics in a classical, point-wise manner while systematically avoiding functorial descriptions. 

The mathematical definitions, constructions, and results reproduced below are standard results in set theory and categorical logic. The purpose of the following is expository and reconstructive: I translate Badiou's presentation into contemporary categorical language rather than claiming new mathematical results.

\begin{dfn} [Category]
A category $\ce$ consists of a collection of objects $\Ob(\ce)$ and, for every pair of objects $A, B \in \Ob(\ce)$, a (not necessarily proper/ZFC) class of arrows $\Hom(A, B)$. It is equipped with associative composition maps $$\Hom(A,B) \times \Hom(B, C) \longrightarrow \Hom(A, C)$$ that connects any two compatible arrows. Finally, each object $X \in \Ob(\ce)$ retains an identity-arrow $\Id_X: X \to X$ with specific identity laws reflecting its compatibility with the composition operator. \end{dfn}

Furthermore, one needs a concept of a relation not just within a single category, but also between them. 
 
\begin{dfn} [Functor] 
A covariant functor $F: \ce_1 \to \ce_2$ between two categories is a mapping $$F: \Ob(\ce_1) \to \Ob(\ce_2)$$ and $$F: \Hom(A, B) \to \Hom(F(A), F(B))$$ that preserves identity and composition. It establishes a relation between two categories of compositions in the same manner that a function relates two different sets of consistencies. 
\end{dfn}

Although for a so-called `small category', a functor could be classically defined as a function between sets of objects and arrows, the structuralist paradigm treats the relations between categories as primary, thereby avoiding the formal dominance of external definitions anchored on point-sets---a new \textit{ethos} Badiou specifically does not fully register. 
This \textit{new}, `Kantian' paradigm particularly attempts to relativize identities in a purely structural manner: it approaches objects---most notably sheaves, stacks, and groupoids---through the matrix of how they map into each other. As I have previously argued\footnote{Veilahti, `Alain Badiou's Mistake', 2013.},
this transformation manifests deep analogues to Kant's transcendental philosophy, unlike Badiou's `Platonic' perspective reluctant to abandon the static ontology of consistencies implicated by set-theoretic formalism.

This reduction, which notably limits how semantics and ultimately the subjectivities inherent to mathematical inquiry are conceived, becomes the most obvious when Badiou attempts to discuss not only intra-worldly relations between objects but relations between worlds themselves\bk ---for example, the geometric morphisms giving rise to the $2$-category $\Top$. Instead of tracing the structural shift from the `despotic' truth-form $\truth$ and the localic extension $T$ to relational constructs like the geometric  subobject classifier $\Omega$ in a general topos, Badiou can understand the meaning of `points' only through the question of true or false---i.e. in a form resembling Kierkegaard's phrase `either/or'---rather than endorsing more radical or geometric points made expressible by category theory: points endowing \textit{internal} endomorphisms invisible at the level of element-wise, incorporeal presentation. 

This relates to Grothendieck's crucial insight, which relied on his observation that if morphisms of varieties were considered via their `adjoint' field of functions, geometric epimorphisms would correspond to monomorphisms of fields. This analogue was not inevitable, however: the category of fields did not possess an operator analogous to pullbacks (fibre products) \textit{unless embedded within a larger category of rings}, in which pullbacks have the co-dual of tensor products $\otimes$. While a covering space is locally `split'---as mathematicians term it---the same was not true for the dual in the category of fields, where every morphism is monic. This led Grothendieck to realize the necessity of replacing classical neighborhoods $U \into X$ with the more general category of arrows $U \to X$ that are not necessarily monic.

Grothendieck was then able to bypass the problem of `split' morphisms by generalising 
the `transcendental indexing' (i.e., coordination) of topological neighborhoods---moving from the 
external Heyting algebra $\Ou(X)$ of a Kuratowski space to a general category $\ce$. 
The synthetic use of \textit{sieve} structures enabled the generalization of 
sheaf-theoretic tools into the category of presheaves $\sets^{\ce^{op}}$, where 
the index category $\ce$---unlike Badiou's localic indexing $T$---may possess non-trivial endomorphisms. 

Topological categories such as Grothendieck topoi retain a truth(-value) object 
$\Omega$ structurally analogous to Badiou's `transcendental' object $T$, but it is not characterised by explicit `values' (as in a set) that 
could be graded as an explicit collection of `degrees'. This is precisely because of the endomorphisms endogenous to the transcendental object $\Omega$, conjugating the hinterland of any specific `truth-value'. Instead, the generalized concept of this object is characterised by the so-called \textit{universal property} of the unique arrow called $$\true: 1 \to \Omega,$$ 
whose pullback property classifies the subobjects of any object (monic) in that 
category. In particular, the endomorphisms of the truth-value object $\Omega$ 
are closely connected to the Grothendieck topology\footnote{Johnstone, \textit{Topos Theory}, 1977, p. xiv.}. 
Grothendieck's core idea was not to force out such endogenous torsion in the first place, like Badiou's formalism does, 
but rather to approach the architecture from the standpoint of only locally, and not generally hierarchical snapshots or `sieves', giving rise to a definition of a sheaf that 
prevents a direct recourse to a global representation of relations as explicit, 
set-theoretic functions.

\begin{dfn} [Grothendieck Topology] \label{sieve}
A \textit{sieve}---in French a `crible'---on an object $C$ is a covering family of $C$ that is downward closed\footnote{I use sieve in the standard order-theoretic form appropriate for the posetal sites considered below. I focus on posetal sites rather than more general categorical settings in order to maintain the connection to Badiou's work, which relies on an external Heyting algebra. I seek to reproduce, in categorical form, the theorems and proofs originally formulated in this posetal setting, which Badiou further restricted to the context of localic sets.}.
 
A Grothendieck topology on a category $\ce$ is a function $J$ assigning a collection $J(C)$ of sieves for every $C \in \ce$ such that:
\begin{enumerate}
\item the maximal sieve $\{f \mid \textrm{for any } f: D \to C\} \in J(C)$,
\item (stability) if $\Si \in J(C)$, then
$h^*(\Si) \in J(D)$ for every morphism $h: D \to C$, and
\item (transitivity) if $\Si \in J(C)$ and $R$ is a sieve on $C$ so that  $h^*(R) \in J(D)$ for all arrows $h: D \to C$, then $R \in J(C)$. 
\end{enumerate}
\end{dfn}

The stability axiom particularly specifies that the intersection of two sieves in a topology also lies within it. Given the latter two axioms, a topology itself emerges as a sheaf on the maximal topology consisting of all sieves. Although the definition of a Grothendieck topology does not generally refer to a basis, it is often useful to operate with one as it makes the idea of what a non-localic topos might look like more tangible.

\begin{rmk}
The notion of a Grothendieck topology enables the generalization of the concept of a \textit{sheaf} to a broader categorical terrain. If $\esp$ is the category of topological spaces and continuous maps, sheaves on $X \in \Ob(\ce)$ are equivalent to the full subcategory of $\esp/(X, \ce)$ consisting of local homeomorphisms $p: E \to X$. The corresponding pair of adjoint functors $\sets^{\ce^{op}} \rightleftarrows \esp/(X, \ce)$ yields an associated sheaf functor $\sets^{\ce^{op}} \to \Sh(X)$, demonstrating that sheaves can be conceptualized either as presheaves satisfying an exactness condition or as spaces with local homeomorphisms into $X$\footnote{Ibid., pp. 10--11.}.
\end{rmk}
 
\begin{rmk}

The category of sheaves $\Sh(\ce,J)$ is a full subcategory
of the presheaf category $\sets^{\ce^{op}}$. These categories
are equivalent when $J$ is the trivial Grothendieck topology,
in which only maximal sieves are covering. In general, however,
the sheaf condition imposes additional constraints on presheaves,
and the two categories need not be equivalent. When no topology
is specified, the category $\sets^{\ce^{op}}$ is understood as
the category of presheaves, rather than of sheaves for an
unspecified Grothendieck topology.
\end{rmk}

\begin{rmk} [Spaces] 
The notion of a sheaf extends the `point-wise' representation of a `classical space'---traditionally conceived as a carrier of consistent functions---into a `diagrammatic' paradigm that composes the functorial relations underpinning such carriage. A sheaf thus `extends' the structure of a space by designating points not as `pre-established' entities or elements. Rather, they are synthetic results emerging from hierarchical, but in Grothendieck's work \textit{contextualized} procedures of \textit{localization}\footnote{Veilahti, `Alain Badiou's Mistake', 2013, p. 17.}.

To exemplify this, suppose I know nothing about an entity $X$, whether it 
constitutes a space or not. If I am told that for any `given' multiple $U$ there 
exist arrows (relations) $X(U) = \{U \to X\}$, I would not yet be able to discern 
the geometric structure of $X$. Such information is required, but I, in addition, 
need to know how different ways of probing $X$ via these test spaces relate to one another. 
That is to say, I must also know how the elements of two different sets of generalized elements, $X(U)$ and 
$X(V)$, remain \textit{compatible} with each other within the intersection $U \cap V$. 
Consequently, for every pair of test spaces $U$ and $V$ that we come up with, and for 
every arrow $\phi: U \to V$ mapping these into each other, we require information on 
how $\phi$ maps the generalized elements of $X(V)$ into generalized elements of $X(U)$ through the composition 
$U \to V \to X$. The sheaf conditions turn out to be enough, if satisfied, to guarantee 
that this structural information amounts to an essentially unique object $X$---an object called a \textit{sheaf}.

In fact, the Yoneda Lemma establishes that sheaves carry all relevant information about the underlying  `space' $A$, which can be uniquely associated with 
its corresponding representable presheaf functor $h_{A}: B \mapsto \Hom(B, A)$. Given 
a sheaf $\F$, the elements $\F(A)$ correspond to natural transformations between the 
functors $h_A$ and $\F$, depicted via the transformation structure 
$$\xymatrix@C+2pc{\ce \rtwocell^{h_A}_{\F} & \sets}.$$ 
This infrastructure culminates in the \textit{Yoneda Embedding}:\footnote{For a full proof 
of the embedding, see Mac Lane \& Moerdijk, \textit{Sheaves in Geometry and Logic}, 1992.}
$$\mathbf{y}: \ce \to \sets^{\ce^{op}},$$
which embeds the index category $\ce$ into its corresponding Grothendieck topos 
$\sets^{\ce^{op}}$, mapping the objects of $\ce$ directly into \textit{representable 
pre-sheaves}.
\end{rmk}

\begin{rmk} [Sheaf of Functions]\label{functions}
While some (pre)sheaves thus `represent' ordinary (topological) spaces $X$, the general paradigm of a sheaf is to designate a \textit{categorically valid} composition of an object that is \textit{analogous} to an entire \textit{set of functions}. For example, one may consider continuous functions on $X$ and the `sheaf' of continuous functions $\mathcal{C}\mathpzc{ont}(U)$ consisting of continuous maps defined on the subspace $U \subset X$. The category-theoretic objective is to arrive at an analogous class of functions without treating $U$ as an \textit{incorporeal subset} of $X$ as in Badiou's formalism, but rather considering any arrows $U \to X$ that inscribe a categorical relation between objects. A sieve is then precisely such a synthetically hierarchical, \textit{contextualized structure of forcing} that allows the sheaf-theoretic \textit{compatibility} conditions to be defined and determined, but only for one context at a time, and not for the entire object simultaneously. 
\end{rmk}

\begin{rmk}[Points]
Because not all sheaves are representable, the Yoneda Embedding demonstrates that \textit{not all points have to be representable} in the traditional, set-theoretic sense of the term. In other words, not every `point-like' sheaf is totalisable by a split decision of the type `either/or'; instead, it may retain an \textit{endogenic life} that remains indiscernible by the set-theoretic means of description. 

Under this view, the `point' itself transforms into an `object' whose consistency is `inaccessible' in the localic setting. If one were instead to conceptualize a point in purely `atomistic' terms---prohibiting any interiority of the point through an absolute decision---one would, in fact, have already breached the Kantian limit by relying on a `Platonically' idealist image. This very image crosses the transcendental boundary by prescribing a restrictive, incorporeal and a global \textit{path} as the exclusive mode through which \textit{atomistic} points are permitted to relate to one another. Category theory, in contrast, depicts points in radically contextual, non-reductive terms. 
\end{rmk}

\section{Points: How Badiou Denies Their Geometry? }
\label{sec_topology}

Badiou's discourse on `points', contra the categorical designation, illustrates these global restraints by relating strictly to the set-theoretically bounded regime of locales. It thus effectuates particular, `pre-established' constraints relevant to Badiou's ontological, `Platonic' definition of points. When he raises the question of \textit{subjectification}, his only recourse is therefore to reproduce the idealist discourse of `being', strictly framed in a two-valued setting, instead of treating them more relationally. To Badiou, a `point of truth' is thus an external complete Heyting algebra morphism $T \to \truth = \Two$, where the latter is the canonical Heyting algebra structure. 

Let me first overview the standard, topos-theoretic discourse on `points' before reframing this in a Badiouian, localic setting. I must first define a \textit{geometric morphism}, which is established as an \textit{adjoint pair} of functors that effectively relate \textit{two different topoi} (i.e. as a map between two `worlds').  

\begin{dfn}[Adjoints]
Functors $f_* : \Es \to \Fs$ and $f^* : \Fs \to \Es$ between topoi are adjoint (denoted $f^* \dashv f_*$) if they satisfy the following universality condition: for every object $X \in \Ob(\Fs)$ and $Y \in \Ob(\Es)$, there exists a canonical bijection 
$$\Hom_{\Es}(f^*(Y), X) \cong \Hom_{\Fs}(Y, f_*(X)),$$
where the left adjoint $f^*$ is required to be left exact.
\end{dfn}

\begin{dfn}[Geometric Morphism]
A geometric morphism between two topoi $\Es$ and $\Fs$ is an adjoint pair $(f_*, f^*)$ consisting of a direct image functor $f_* : \Es \to \Fs$ and a left-exact inverse image functor $f^* : \Fs \to \Es$. 
\end{dfn}

\begin{dfn}[Point in a Topos]
Given a $\sets$-topos $\gamma : \Es \to \sets$, a point is a geometric morphism $p : \sets \to \Es$. It is not necessary that any points exist---this is the case in what \LW\ refers to as an `atonic world'---nor that they are sufficient to determine the structural identity of objects within the topos. For a class of points $K$, by contrast, one says that $K$ retains \textit{enough points} if the family of functors $p^*$ determines isomorphisms; that is, in the case that an arrow $f : X \to Y$ is an isomorphism if and only if all set-theoretic pullbacks $p^*(f)$ are bijections\footnote{Johnstone, \textit{Topos Theory}, 1977, p. 224.}.
\end{dfn}

\begin{rmk}
A point in topos theory is thus a concept that relates Badiou's idealist topos $\sets$ with the general, diagrammatic expression of a categorical topos $\Es$. As a morphism between two worlds, it is a translation of endogenic truth into a situation that is obstinately binary. 

This relation is directionally opposite to the one expressing the question of \textit{materialization}, in the sense of my previous article, as the latter is governed by the global sections morphism $\gamma$, determining how precisely a topos is `defined' over $\sets$. However, even if it bridges set-theoretic `ontology' with categorically defined `topology', the morphism $\sets \to \Es$ does not need to relate to logical descriptions in any intelligible manner. Only in the case of locales, i.e. Badiou's worlds, there is a direct connection to the `degrees' of truth, which determines the point directly as a set-theoretically structured \textit{choice} (or a `point of truth' as Badiou terms it). 
\end{rmk}

\begin{dfn}[Morphism of Locales]
If the general definition of a point takes it as a geometric morphism $\sets \to \Es$, 
one needs to justify how it transforms into an algebraic question on Heyting algebras. 
Let $T_1$ and $T_2$ be two extensive complete Heyting algebras of extensive degrees, 
such as lattices of numbers or other ontologically `familiar' objects. If one considers 
a frame homomorphism $\phi : T_2 \to T_1$ preserving the complete Heyting algebra structure, it amounts 
to a map of sheaves---in other words, a family of natural transformations inducing 
a geometric morphism
$$\sets^{T_1^{op}} \to \sets^{T_2^{op}}$$ 
by mapping a sheaf $F : T_1^{op} \to \sets$ to the functor $F \circ \bar \phi : T_2^{op} \to \sets$. 
Here, $\bar \phi$ is the induced map which, for $q \leq p$ in $T_2$ (corresponding to the 
arrow $p \to q$ in the opposite category), maps $\bar \phi(p \to q) = (\phi(p) \to \phi(q))$\footnote{In 
fact, given a functor $\ce \to \D$ between any categories, not necessarily Heyting, the 
functor induces a geometric morphism $f : \sets^{\ce^{op}} \to \sets^{\D^{op}}$ for which 
$f^*(\F) = \F \circ f^{op}$. See Mac Lane \& Moerdijk, \textit{Sheaves in Geometry and Logic}, 
1992, p. 359.}. In the opposite direction, such a morphism is easily seen to entail that 
the complete Heyting algebra structures are preserved. Thus, a point in the case of a 
locale may be defined as a bivalued $\wedge$-$\Sigma$-map: 
$$T \to T_{\sets} = \Two \quad (=\truth).$$
\end{dfn}
Here, the operations on the truth-value object $\truth$ are defined as $$\textrm{true } \wedge \textrm{ true } = \textrm{ true},$$ $$\textrm{ true } \wedge \textrm{ false } = \textrm{ false}$$ and $\Sigma B' = \textrm{ true}$ unless $B' = \{\textrm{false}\}$, meaning it consists of the minimal 
degree alone.

In the case of a localic point $\sets \to \Es$, where $\Es$ is in fact the world of $T$-sets, the point-morphism thus reduces to a map 
$$\phi : T \to \Two$$ which satisfies $\phi(p \wedge q) = \phi(p) \wedge \phi(q)$ and 
$\phi(\Sigma B) = \Sigma \{\phi(p) \mid p \in B\}$. 
Badiou, in effect, masks the `sheafified' insight of categorical geometry. Rather 
than following Grothendieck's locally hierarchical `snapshots', Badiou aims to 
re-subjectify a globalist, absolute concept of points. In other words, whereas a 
general Grothendieck topos allows categorical `torsion'---which can only be filtered context-wise 
via synthetic sieve structures to enable analytic forms of mathematical determination---in the case of \textit{locales}, such filtered structures can be extended to the world as a whole. Only then a binary outcome of a point can be determined---not just against a particular localising sieve, 
but for the entire situation. 

Expressing this differently, to Badiou\footnote{Badiou, \textit{Being and Event}, 2006, 
p. 401.} the question of `topology' reduces to the question of a (local) place whose 
determination can be executed globally. This follows the despotic doctrine of the strong, 
the atomistic `postulate of materialism', which rhetorically inverts the deductive logic 
of the statement `where there's choice, there's a place'---the topology of this place 
is structured as a \textit{despotic point} of truth. Only the rigid choice `either/or' 
allows `points' to emerge in the Badiouian sense, whereas Grothendieck's genius was precisely to escape such a \textit{split} question of points and neighborhoods in the first place. 

\begin{rmk}[Despotic Points] 
These logical points of truth, as formalized in \LW, are termed despotic because they are \textit{globally} 
determined by the set-theoretic truth-form $\truth$, rather than operating  
within locally effective, Grothen\-dieckian sieve structures.
\end{rmk}

The radical shift of the categorical expression of points, and where the Yoneda Lemma serves as the primary engine, was thus to understand why the set-theoretic totalization of points is unnecessary. Even through a point is seemingly isolated or 0-dimensional, it may capture 
 internal morphisms. A point thus needs not be treated as an in-itself absolute object, whose intrinsic contents the Kantian paradigm refuses to consider, but as something dynamically situated within its own diagrammatic compositions. 

This geometric, and in many ways intensive understanding of points is crucial also to the topologisation of Cohen's procedure as I shall illustrate: the primary role of points is not to constitute a topological substance, but they too operate procedurally, reflecting localizing and enforcing acts over a topology. 
By contrast, Badiou's reactionary point, which is a localic reaction that forces an annihilation of this geometric intensity, serves an opposite function. It is not mathematically illegitimate, but a specialized instance that becomes the engine of his philosophical limitations. 

This partialism is most obvious in how Badiou believes that being itself, and thus also topology, is only accessible through its parts, rather than relationally. 
Thus, again, he reads the categorical idea of topology inside-out: first he denies the internal geometry of points, and then recollects their split, despotic representations as a set on which he imposes a classical Kuratowski topology, only to hide his uncertainty by then trying to cursorily reconnect with diagrammatic, functorial compositions.

\begin{rmk}[Badiou's Reterritorializing Topology of Points] 
To (re)con\-struct, indeed, an `ontological' or `singular'\footnote{The term singular in this occasion is used following the standard mathematical terminology, referring to Kuratowski topology, as opposed to the philosophical concept of the singular used elsewhere in this article, referring to the incongruence of any mundane representation. Yet these two meanings, at least in Badiou's setting, are indirectly related as the singular point-set topology specifically avoids a more relational description, while solely retreating to an `ontological', elementwise definition.} topology anchored on the representation of the set of points as a point-set, Badiou\footnote{Badiou, \textit{Logics of Worlds}, 2009, pp. 441--442.} writes $P_p$, termed a \textit{positivation}, for the set of points 
$$\pi(T) = \{\phi: T \to \mathbf{2} \mid \phi \textrm{ is a surjective } \cap\text{-}\Sigma\text{-function}\}\footnote{Badiou idiosyncratically employs the set-theoretic intersection symbol $\cap$ alongside the supremum operator $\Sigma$ to denote the operations of the transcendental lattice. In standard topos-theoretic and categorical frameworks, this corresponds to a complete Heyting algebra morphism preserving the finite infima $\wedge$ and arbitrary suprema $\bigvee$. I retain Badiou's notation here to underscore his structural reliance on point-wise, set-theoretic analytics.}$$ 
and for each $p \in T$ a `positivation of $p$' as a subset of $\pi(T)$, denoted by 
$$P_p = \{\phi \in \pi(T) \mid \phi(p) = 1\}.$$ 
Thus, like Badiou\footnote{Ibid., p. 418.}  observes, `the `active' parts of a group of points are constituted by the positivations' and `we will call the \textit{interior} of a group of points the active part of the group, namely the union of all the positivations which are parts (or subgroups) of the group in question'. Formally, this means that for a subset $Y \subset \pi(T)$, one defines $\Tint Y = \bigcup_{P_p \subset Y} P_p$. These so-called positivations provide a basis for defining a topology on $\pi(T)$. 
\end{rmk}
 
 In later sections I will explore what becomes visible once this specialization that conceals the internal geometry of points is no longer treated as a necessity. Yet, in Badiou's work, their recollection as a set serves as an ontological anchor---the set of points becomes a point-set that, in essence, exiles the question of internal structure from the points themselves. In fact, the reterritorialization of points as a set serves precisely as a way of capturing a topos itself `ontologically'---even if, outside Badiou's partialist framework, this image $\pi(\Es)$ of a topos captured through its points is unfaithful to the structure of the topos itself.

\begin{rmk}\label{closure_remark}
The claim about Badiou's partial understanding of topology is further supported also by how his own definition  deviates from standard mathematical descriptions. This is evident as he attempts to reconnect with a more functorial description, perhaps to claim that it should not bring about anything novel, which though is obviously true only in his localic setting. 
Indeed, to `equip oneself with a more topological intuition and to treat the degrees of the transcendental as operators of the localization of multiple-beings', Badiou thereby elaborates with an alternative route: a topology defined by an  \textit{interior}-operator $\Tint: \Power \pi(T) \to \Power \pi(T)$. 
In formal terms, for $\Tint$ to yield a Kuratowski topology\footnote{In a traditional Kuratowski topology, a topological space $X$ is defined as a set of points itself together with a collection of open subsets closed under arbitrary unions and finite intersections, containing both the empty set and $X$ as a whole. An interior operator maps any subset $Y \subset X$ into the union $\Tint Y = \bigcup_{U \subset Y \textrm{ open}} U.$}, it is required to satisfy, for all $A, B \subset \pi(T)$, the  conditions: $\Tint(A) \subset A$, $\Tint(A) \cap \Tint(B) = \Tint(A \cap B)$, and $\Tint(\Tint(A)) = \Tint(A)$. Therefore, even if advancing a more functorial intuition, Badiou's description never quite escapes the traditional Kuratowski site defined over explicit, `ontologically' determined sets. 

In effect, Badiou's `functorial' intuition based on $\Tint$, but still anchored to point-sets, appears ultimately misleading. To provide further evidence of this philosophical partialism---even if here only circumstantial, based on his formal confusion---categorical topology is standardly defined as the \textit{closure} operator $j$ that is compatible with arbitrary (read non-finite) intersections, like the double negation $\neg\neg$ operating on $T$ or $\Omega$. By contrast, Badiou's interior operator $\Tint$ is compatible with arbitrary unions: a difference Badiou symptomatically omits\footnote{Obviously, the operator $\Tint$ itself easily translates into a Kuratowski closure: $\bar A = \pi(T) \smallsetminus \Tint(\pi(T) \smallsetminus A),$
but Badiou makes some formal mistakes as he demonstrates its compatibility with unions and intersections, indicating deeper confusion over the significance of closure as a topological engine.}.

\end{rmk}

\section{Forcing the `Generic'---Dominations and Sieves Revisited}

Given that category theory engenders a more geometric concept of points, which play a central part in \LW, I will explicitly bridge this understanding with the theory of subjectivity presented in \BE. In this work, Badiou\footnote{Badiou, \BE, 2006, pp. 420--423.} demonstrates Cohen's proof of the `errancy of excess'---the independence of the continuum hypothesis: 
\begin{hyp}[Georg Cantor, 1878]
There is no set whose cardinality is strictly between that of the integers and that of the real numbers.
\end{hyp}

 In short, Badiou argues both for and against the claim that in the subject-language of a \textit{specific model} of set theory, the cardinality $|\Power(\omega_0)|$ of natural numbers can `surpass an absolutely indeterminate cardinal $\delta$ given in advance'. The cardinality of such an object $\delta$ is obviously interpreted in terms of the language pertaining to the `inhabitant of the situation'---\textit{qua} the particular interpretation function---since `ontologically' speaking such a claim would render the situation inconsistent and thus paradoxical. In the proof that follows, I employ the (internal) cardinality $\delta = |\Two^N| = \Power_{int} N$ for the so-called internal natural number object $N$, but the construction is fundamentally the same for any $\delta$. In set-theoretic terms, following Badiou's set-theoretic description, the procedure involves the consideration of the so-called `set of conditions' $\copyright$.  

\begin{dfn}[Condition] 
A condition is a certain finite subset of the set of all triplets $\langle a, n, \beta\rangle$ where $\beta \in \{0, 1\}$ is a binary number. A condition is specifically a subset where for any given $a \in \delta$ and $n \in N$, both $\langle a, n, 0\rangle$ and  $\langle a, n, 1\rangle$ cannot simultaneously belong to the set of conditions. 
\end{dfn}

\begin{rmk}
It is easy to see how these triplets, or (finite) partial functions of the form $\delta \times N \to \mathbf{2}$, can be used to engender new `points of truth' in the sense of \LW. Namely, since only one of the two triplets $\langle a, n, 0\rangle$ and $\langle a, n, 1\rangle$ may be contained within a single condition, each pair $\langle a, n\rangle$ corresponds to a single binary value. As Badiou\footnote{Ibid., p. 420.} phrases, the final index, either 0 or 1, signifies whether `$a$ is \textit{paired} with $n$'.
\end{rmk}

\begin{rmk}
Combined together, the set of conditions $\Cond$, or $\copyright$ as Badiou terms it, forms a poset.
Indeed, each condition $p \in \cond$ is a finite subset of $\delta \times N \times \Two$, and so we can define an order-relation as: 
$$q \leq p \textrm{ if and only if } q \supset p.$$
While existing as partial functions, these conditions $p$ can be further conceptualized as local topological domains over which a desired map used in Cohen's proof, $\gamma: \delta \times N \to \Two$, is approximated as a sheaf $F$. A set of local sections $F(p) = \{q \in \Cond \mid q \leq p\}$ would then consist of local sections so that a stronger condition $q \leq p$ would define a smaller and more precise window of approximation of $\gamma$, defined through a `restriction map' $F(q) \to F(p)$ (which is actually an inclusion). Moreover, one can define a compatibility-condition $p \dagger q$ consistent with the formalism of \LW: $p$ and $q$ are compatible if they  approximate a shared partial function $\gamma$, that is, when they admit a common stronger condition $r \leq p, q$, or equivalently when $F(p) \cap F(q) \neq \emp$\footnote{Badiou idiosyncratically employs the dagger symbol $\dagger$ to denote the structural compatibility of forcing conditions. In standard set-theoretic and topos-theoretic frameworks, this relation is conventionally written as $p \parallel q$, signifying that the conditions possess a common lower bound $r \leq p, q$ in the poset. I retain Badiou's notation here to preserve textual fidelity to \textit{Logics of Worlds} while mapping its intra-worldly consequences.}.

\end{rmk}

The poset $\cond$ is not directly a Heyting algebra, or `transcendental indexing' in the terminology of \LW: given that all conditions are finite sets, there is no `minimal degree' $\mu$ compatible with the order-relation. Instead, a Heyting structure emerges from considering downward closed subsets of $\cond$---the so-called `sieves'. Badiou himself engages apparently similar constructs, which he calls `dominations', but in reality they are not necessarily downward closed nor defined relative to specific conditions, illustrating his formal imprecision: 

\begin{rmk}[Badiou's Dominations]
In \textit{Being and Event}, Badiou\footnote{Ibid., p. 369.} defines `dominations' as subsets of the forcing poset $\copyright = \Cond$ satisfying the rule: 
$$ \text{if } p \notin D, \text{ there exists a } d \in D \text{ such that } p \subset d. $$
Badiou further calls them `dense' if they are downward-closed, which in standard vocabulary instead means that they are sieves---in fact, as Badiou defines the dominations, they are always dense in the proper mathematical vocabulary, where density is formalized in the sense of the $\neg\neg$-topology of the associated Cohen topos (see the next section).  
\end{rmk}

\begin{dfn}[Sieves]\label{sieve_cover}
The standard definition utilized within Cohen's actual proof indeed defines these domination-structures as \textit{downward-closed} families---a property Badiou conflates with density---meaning that for any $d' \leq d \in D$, it strictly follows that $d' \in D$. This property renders them  \textit{sieves} in sheaf theory, an essential attribute that can be integrated into Badiou's exposition without structural alterations to the underlying forcing argument. 
Crucially, the dominations in the actual proof that Cohen conducts---to preserve the modal transformations between the ground model and the generic extension---need to be defined \textit{relative to a specified condition} $p$, but not necessarily $p = 0$ as in \BE, resulting in the localized density condition:
$$ r \leq p \quad \Longrightarrow \quad \exists d \in D: \quad d \leq r. $$
Such constraints remain necessary for Cohen's transfinite framework, reflecting Badiou's formal imprecision when employing a less accurate designation of dominations. A sieve $D$ satisfying this localized density condition relative to an initial condition $p$ is formally said to \textit{cover} the condition $p$---a crucial property to defining the closure operator in the structural topology of Cohen's argument. 
\end{dfn}

\begin{rmk}[Forcing Points]
Cohen's original idea is to consider subsets $$\gamma(n) = \{ \langle a, n, 1\rangle \in \female \} $$ for a suitable `generic' sieve $\female \subset \cond$---a subset Badiou\footnote{Badiou, \BE, 2006, p. 335.} calls `indiscernible' because it is inaccessible `within the situation' of the original model of ontology. Though expressing the argument in a philosophically confusing way, while failing to explicate that the `inhabitant' itself changes alongside the model, the idea nonetheless bridges with $\LW$ in the sense that the generic sieve $\female$ directly corresponds to a `point of truth' by emerging as the exponential transpose of the map $\female: \cond \to \Two$. It `groups together all the terms of the situation which are positively connected to the event', or rather, the ontological change of the situation. A new binary point of truth thus effectively appears---or rather becomes expressible---precisely as $\cond$ is \textit{split} into two parts by $\female$. In other words, `[t]he essence of the demonstration consists in revealing dominations in $S$, which must consequently be `cut' by the generic part $\female$'\footnote{Ibid., p. 423.}: obviously a process conceivable only on the assumption that being itself is first conceived in binary terms. 
\end{rmk}

Now, this idea is far from stating that the despotic, split truth point of truth would change the status of a given, ontological situation or model. Instead, the splitting, procedurally indebted to the axiom of choice, is strategically used to \textit{force the expression} of new points which in this sense are mathematically \textit{effective} rather than causal. And, given their expressive nature, the process itself---unlike the specific outcome $\female$---resonates with Grothendieck's insight that was precisely to avoid assuming the `ontologically' split presentation of algebraic varieties, but instead to only approach such internally automorphic sites through \textit{local hierarchies}, or sieves, effectuated on specific contexts.

The splitting, in effect, should be first and foremost understood procedurally, rather than atomistically. 
It is precisely the universal property of the exponential transpose---the transcendental gaze  Badiou struggles to engage in \textit{Logics of Worlds}---that effectuates a procedural ramification around the map central to Cohen's demonstration: 
$$\gamma: \delta \times N \to \mathbf{2}. $$ While this relation bears a formal resemblance to the structure analyzed by Badiou\footnote{Ibid., p. 420.}, a crucial difference emerges: Badiou's construction operates via the exponential transport $\gamma(n)$\footnote{Ibid., p. 422.}, whereas my proof works in the opposite direction by initiating the procedure from $\delta$. Consequently, I employ the opposite transpose $\gamma: \delta \to \Power N$, yielding the specification:
$$\gamma(b) = \{n \mid \gamma(b, n) = 1\}.$$ 

The supposition that $\gamma$ constitutes a monomorphism entails that if $b \neq b'$, then $\gamma(b, n) \neq \gamma(b', n)$ for some $n$. In the base model $\Ss$, by Cantor's classical diagonal argument, no such monic arrow can exist. However, much like the infinite knowledge accessible from the `generic' position in \textit{Being and Event}---an architecture unknowable to a finite inhabitant and accessible only to an ontological alien---this field can be filtered by a finite series of compatible, ordered elements; that is, finite approximations.

Badiou, in contrast, views the structures of forcing from a globalist point of view, as if causally regulating the situation as a whole---the precise source of event-related anguish. This becomes his `fundamental law of the subject', which, perhaps, should rather be renamed as the idealist or `Platonic' law forcing the splitting on the dialectical whole: 
\begin{quote} ``If a statement of the subject-language is such that it will have been veridical for a situation in which a truth has occurred, this is because \textit{a} term of the situation exists which both belongs to that truth (belongs to the generic part which \textit{is} that truth) and maintains a particular relation with the names at stake in the statement. This relation is determined by the encyclopaedic determinants of the situation (of knowledge).''\footnote{Ibid., p. 401.}
\end{quote}

\section{Presenting the `Errancy of Excess' Categorically} 

To turn now to the actual demonstration, a critical fact Badiou never states---and that reveals its connection to the theory of points---is that the `generic filter' or sieve $\female$ is nothing other than the exponential transpose of $$\gamma: \delta \times N \to \mathbf{2}, $$ thus directly compatible with the notion of points in \LW. Cohen's procedure thus reduces to a matter of demonstrating that within a specific, extended model $S(\female)$ of set theory that is deliberately constructed to \textit{enforce} the emergence of the point $\female$,  one can---depending on the construction\bk ---\bk either prove or disprove the continuum hypothesis: that `\textit{there are at least as many parts of $\omega_0$ of the type $\female(n)$ as there are elements in the cardinal $\delta$}'. Again, this statement pivotal to Cohen's proof is not approached as an `ontological' claim---though Badiou makes it appear so---but as a statement expressible \textit{within} the modelled language  of $S(\female)$\footnote{Ibid., p. 420.}. 

Given the denumerability of the model $S$ that Cohen's proof assumes by definition, `ontologically' $\delta$ is clearly denumerable as well. But as \textit{interpreted} by `an inhabitant of $S$', $\delta$ would appear non-denumerable should there not exist a bijection between $N$ and $\delta$ expressible in the internal language $\Expr_S$ interpreted within the situation through the so-called interpretation function $\Int_S$. Given this dual status of language, which effectuates a `paradoxical' experience internal to the situation, the result is non-causal with respect to `being' as such, though linguistically effective. Thus, the leap to topos theory is rather natural, where the `internal' model retains a less idiosyncratic or superimposed character, deliberately avoiding that implicit confusion on which Badiou's somewhat questionable philosophical conclusions about the `errancy of excess' rely. 
 
 Next, I will overview the necessary categorical ingredients to present an analogous argument in the category-theoretic framework. First, the `ontological' model $S$ is replaced by a base topos $\Ss$ that possesses all the relevant properties of the binary, set-theoretic world $\sets$---a topos which satisfies a categorical version of the axiom of choice (AC), fulfils the structurally generated (SG) condition, and has a natural number object. One may conceptualize $\Ss = \sets$.

Obviously, the continuum hypothesis requires a well-defined notion of a denumerable infinity, for otherwise the hypothesis itself cannot be formally stated. This is precisely the objective of designating a natural number object $N$---not as an ontological construct, but specifically as a categorical, relational expression. 

\begin{dfn}[Natural Number Object] \label{naturalnumberobject}
A definition due to Lawvere defines a \textit{natural number object} in a topos $\Es$ as an object $N$ (reflecting the arithmetic space of natural numbers) equipped with morphisms $$\xymatrix{1 \ar[r]^o & N \ar[r]^s & N}$$ such that for any diagram $\xymatrix{1 \ar[r]^x & X \ar[r]^u & X}$ in $\Es$, there exists a unique morphism $\xymatrix{N \ar[r]^f & X}$ making the following diagram commute:
$$\xymatrix{
1 \ar[r]^o \ar[dr]_{x} & N \ar[d]^f \ar[r]^s & N \ar[d]^f \\
& X \ar[r]^u & X
}$$
The arrow $s$ corresponds to the classical succession arrow of $\sets$, mapping $1 \to 2$, $2 \to 3$, 
and so forth---with the crucial distinction 
that $N$ itself does not necessarily consist of global, 
discrete elements. Such a number object is central for the emergence of 
internal algebraic theories endogenic to a topos (see Section \ref{mitchellbenabou}).

In fact, the existence of the natural number object is also philosophically tempting because it reflects a sense of algebraic \textit{finality}: it can be expressed equivalently as the existence of a free $\mathbf{end}$-model functor, where $\mathbf{end}$ denotes the algebraic theory generated by a single unary operation $e$ without equations\footnote{See Johnstone, \textit{Topos Theory}, 1977, p. 221.}. 
\end{dfn}

Provided that the natural number object $N$ is given, one way to express Cohen's procedure would be just a direct translation through suitable geometric morphisms $\Es \to \Ss$. In other words, it is possible to reframe the notion of conditionality so that valid statements internal to $\Ss$ can be interpreted as suitable, \textit{conditionally valid} propositions expressed in $\Es$: for instance, a claim that, for an existing object $X$, there are no epimorphisms from $N \to X$ or from $X \to \Omega^N$ thus falsifying the continuum hypothesis. This approach would preserve the original, set-theoretic presentation of Cohen's procedure, even if through its `internal' interpretation---a strategy that becomes meaningful after the introduction of Mitchell--Bénabou language (see Section \ref{mitchellbenabou}). In particular, it overcomes the need to express alternative models $\Ss'$ from the transcendental position as if they existed independently of the relational reality of $\Ss$.

However, to illustrate the power of the categorical mode of expression, and how it shifts the question of subjectivity away from specific points of truth by operating instead at the level of functorial associations, I will overview the demonstration in a more `transcendental', though relationally structured setting in the next section. To this end, a suitable categorical definition of topology is required. Crucially, it deviates from the classical, elementwise description of the Kuratowski-frame like the one by which Badiou himself seeks to reterritorialize the set of points as reviewed in Section \ref{sec_topology}, essentially reterritorializing also the topos itself whose set of points he seeks to topologize.  Cohen's proof, instead, does not operate on points as such, but on the precursive structure, that is, the poset of `conditions' from which points are forced to emerge only in retrospect, as in Cohen's original procedure, making the geometric aspects of the original proof more apparent.

To be precise, the category $\cond$ is used for the `transcendental', local indexing of sheaves. To this end, it must be equipped with a relationally defined \textit{Lawvere--Tierney topology}: an arrow operating on the internal subobject classifier $\Omega$, which plays the same structural role as the `transcendental indexing' $T$ in Badiou's framework. Formally, such a topology refers to a closure-endomorphism $j: \Omega \to \Omega$ satisfying the structural requirements\footnote{For a Grothendieck topos $\sets^{\ce^{op}}$ where $\Omega(U)$ consists of sieves, the closure operator is defined as a sheaf transformation of the form $j_C(\Si) = \{g \mid g^*(\Si) \in J(\dom g)\}$; that is, the family of maps `covered' by the sieve $\Si$.} (also see Remark \ref{closure_remark}): $$j \circ j = j, \quad j \circ \true = \true, \quad \textrm{and} \quad j \circ \wedge = \wedge \circ (j \times j).$$ 
The specific topology used in the categorical version of Cohen's proof is the so-called double negation topology $j = \neg\neg$.

However, complicating the functorial definition of topology, $\Cond$ itself is not a Heyting algebra or `transcendental indexing' in the sense of \LW, and does not directly retain logical operators; that is precisely why sieves, or what Badiou calls `dominations', are needed to induce a proper topology---the very secret that ultimately bridges the subject philosophy of \BE\ with the underlying topological structures implicit in \LW. 
For the purpose of Cohen's proof, this transcendental object $\Omega$ should be a presheaf over the poset of conditions $\Cond$, and one can simply define its sections over a given condition $p$ as:
$$ \Omega(p) = \{D \mid D \textrm{ is a sieve on } p\}. $$
Likewise, the topology $j$ is a natural transformation that, for each condition $p$, expands an initial sieve $D$ so as to consist of stronger conditions $q$ that $D$ covers (see also Definition \ref{sieve_cover}):  
$$ j_p(D) = \{q \leq p \mid D \textrm{ covers } q\}. $$

\begin{rmk}
Generally, the $j$-sheaves in a topos $\Es$ with topology $j$ are those objects $E$ for which every dense monomorphism $\sigma: A \hookrightarrow E$ induces a canonical bijection 
$$ \Hom_{\Es}(E, F) \cong \Hom_{\Es}(A, F) $$
for any sheaf $F$. This reflects the fact that a $j$-dense subobject completely determines the sheaf under the specific topology, demonstrating the importance of `generic filters' in Cohen's argument\footnote{See Mac Lane, Saunders and Ieke Moerdijk, \textit{Sheaves in Geometry and Logic: A First Introduction to Topos Theory}, 1992. p. 223.}.
\end{rmk}

\begin{rmk}[The Forcing of Classical Worlds]
While the operator $\neg\neg$ is generally defined only internally to the subobject classifier $\Omega$, when the underlying topos constitutes a locale---being both Grothendieck and satisfying the axiom (SG) like Badiou's postulate of materialism demands---the terminal global sections morphism $\gamma: \Es \to \sets$ is bounded. Consequently, the $\neg\neg$-topology translates into an operator acting on the external frame $T = \gamma_*(\Omega)$, bridging directly with Badiou's typology of localic worlds. In what he terms a `classical', discrete world where the law of the excluded middle holds, this mapping reduces to the identity---something achieved in the closure $T_{\neg\neg} = \neg\neg(T)$ itself, which in Badiou's vocabulary is always classical. In the same way, the $\neg\neg$-sheaves form a substructure of all sheaves behaving in a classical, Boolean manner, which is pivotal to understanding the role of topological closure in the forcing of binary truths as Cohen's proof entails.  
\end{rmk}

Although geometric morphisms of the form $\Es \to \Ss$ are generally not logical---otherwise it would be meaningless to attempt to construct two contradictory points of truth like the `errancy of excess' demands---Lawvere and Tierney's construction proceeds by building new models of set-theory: each poset $\cond$ gives rise to a \textbf{Cohen topos} $\Sh(\cond, \neg\neg)$, where `new truths' articulated by its internal logic can be transformed into conditional statements in the original situation $\Ss$. This topos is both Boolean and satisfies the (SG) condition, thereby representing another model of set theory that either negates or confirms the continuum hypothesis. Crucially, the inverse image component of the geometric morphism $\Ss \to \Sh(\cond, \neg\neg)$ fully preserves the structural composition of the natural number object\footnote{Since $g^* : \Es \to \Fs$ is left exact by definition, it preserves the exact sequence defining the natural number object (cf. Mac Lane \& Moerdijk, \textit{Sheaves in Geometry and Logic}, 1992, p. 269). Most crucially, the category of sheaves has a natural number object: if $1$ is the constant sheaf in $\Sh(\ce, j)$, it is isomorphic to the countable coproduct $\coprod_{n \in N} 1$ of copies of the terminal object.}.

\begin{rmk}
Although $\Sh(\cond, \neg\neg)$ is not strictly two-valued---as it fails to yield a despotic truth-object with only two global elements---it can be \textit{forced} into a two-valued model of set theory $\Ss'$ by applying the ultra-power construction analogous to a transfinite recursion, then operating on the filter quotients of the Cohen topos. This yields a variable outcome for the continuum hypothesis within $\Ss'$ depending on the poset $\cond$ used, establishing independence at an `ontological' level of the different entities $\Ss'$. 

With this two-valued Boolean topos $\Ss'$, the global elements of the subobject classifiers $\Omega_{\Ss'}$ directly correspond to the external two-element Boolean algebra $\{\textrm{true, false}\} = \truth$, which is the simplest possible Heyting algebra. Indeed, in every two-valued Boolean topos $\Bs$ (i.e. where $\neg \neg = \id_{\Omega_{\Bs}}$), the maps $\true: 1 \to \Omega_{\Bs}$ and $\neg \circ \true : 1 \to \Omega_{\Bs}$ induce an isomorphism $$1_{\Bs} + 1_{\Bs} \cong \Omega_{\Bs}.$$ 
\end{rmk}

\begin{rmk}
As a minor qualification, even though the quotient construction of the Cohen topos $\Sh(\cond, \neg\neg)$ does not exactly model ZFC, it forms a robust universe of sets in its own right, as Johnstone\footnote{Johnstone, \textit{Topos Theory}, 1977.} observes. The exact ZFC model can be obtained by an argument mimicking the cumulative hierarchy $V_{\alpha}$ for ordinal numbers $\alpha$ internal to the topos, and subsequently executing the same proof on this extended framework\footnote{See Mac Lane \& Moerdijk, \textit{Sheaves in Geometry and Logic}, 1992, p. 279.}.
\end{rmk}

\section{Independence of the Continuum Hypothesis---Revised Proof}

To complete my formal bridge between \textit{Being and Event} and \textit{Logics of Worlds}, the objective of this section is to sketch out the actual proof of the independence of the continuum hypothesis by utilising a rigorous, topos-theoretic framework. The following is not a self-contained reconstruction of every technical step, but it isolates the principal categorical constructions and indicates how the forcing argument is expressed internally to a sheaf topos, illustrating the philosophical relevance of categorical techniques to Badiou's original argument. 

It follows from a standard argument that the category of $\neg\neg$-sheaves within the presheaf topos $\sets^{\cond^{op}}$ is a Boolean topos\footnote{Ibid., p. 272.}. For any subobject $A \into E$, the closure-operator---which extends the closure-operator $j$ previously defined on the sheaf of sieves $\Omega$---is evaluated pointwise as:
\begin{equation}\label{doublenegationsheaf}
(\neg \neg A)(C) = \{x \in E(C) \mid \forall f: B \to C, \, \exists g: D \to B; x \circ f \circ g \in A(D)\}.
\end{equation} 
Because $\cond$ is a poset, the topos $\Sh(\cond, \neg \neg)$ satisfies the internal version of the axiom of choice (AC)\footnote{Ibid., p. 277.}. This critical observation ensures that the independence argument can be faithfully translated back into a proper model of set theory. 
It follows from the fact that $\Sh(\cond, \neg\neg)$ is generated by the subobjects of the terminal object $1$ (SG). Crucially, however, the proof habitually relies on the external, ontological version of the axiom of choice---specifically via Zorn's Lemma applied to the poset $\Sub(I)$---utilized in its practical execution rather than its formal syntax. In other words, the axiom is not established internally to this argument, but rather emerges as an expression of the ontological indebtedness of this sheaf-theoretic construction of the Cohen topos. 

\begin{rmk} 
As the following demonstration underscores, even if category theory expresses an altered transphenomenal domain of mathematics, this same indebtedness characterized by categorical \textit{proofs} is visible in how $\Hom$-sets are designated as classical sets. However, this is not so much a property of categorical expressions themselves, but of their translation into the classical, `Platonic' domain. Whenever mathematicians establish theorems or construct points of truth, 
they remain dependent on ontological arguments---that is, the level 
of binary designation---even if categorical expression as such does not demand or explain this drive for `ontological' proofs, or the need to make a point about Cohen's construct. By contrast, the possibility of expressing Cohen's argument as a relational effect while working solely at the level of internal constructs, like I previously suggested, would reinterpret Cohen's work in a less ontologically predisposed fashion, though this would remain less consequential with respect to set-theory itself and would not indicate the precise procedural vulnerabilities in Cohen's proof that are the source of logical aporia. 
\end{rmk}

To make a point about the `errancy' of set-theoretical models themselves, one shall first prove the following lemma concerning the \textit{internal} version of the Axiom of Choice (AC) within the topos. Crucially, because this demonstration relies explicitly on Zorn's Lemma---which structures an external variant of the Axiom of Choice---the internal satisfaction of the Axiom of Choice within this topos-theoretic model remains dependent on the underlying ontological, designating layer of the proof.

\begin{lem}
A Grothendieck topos generated by the subobjects of the terminal object and for which the class of subobjects $\Sub(E)$ of any object $E$ is a complete Boolean algebra, satisfies the internal version of the Axiom of Choice stating that every epimorphism $X\to I$ splits (i.e., there is a compatible reverse arrow so that $I \to X \to I$ is an identity). 
\end{lem}

\begin{proof}

The crucial part is to demonstrate that for an epimorphism $X \to I$, there is a maximal subobject $M \hookrightarrow I$ such that a section $M \to X$ exists. This requires Zorn's lemma, i.e. the `ontological' equivalent of the axiom of choice typically postulated in set theory. One first considers the set of pairs $(M, s)$ where the original inclusion factors as $M \to X \to I$. This set of partial sections can be ordered by the relation $(M_1, s_1) \leq (M_2, s_2)$ if there is $\iota: M_1 \hookrightarrow M_2$ compatible with partial sections $s_1$ and $s_2$. While in general some of these partial sections might be incompatible, when restricting to a totally ordered \textit{chain} $\{(M_i, s_i) \}$ of partial sections, the completeness of $\Sub(I)$ ensures that there exists a pair 
$(\bigvee_i M_i, s_{\vee})$, which is the upper bound of that specific chain. 

Zorn's lemma, in turn, guarantees that because each chain has an upper bound, there is a maximal element in the entire poset $\{(M, s)\}$. Otherwise, one would always find a larger partial section $$(M_1, s_1) \leq (M_2, s_2) \leq (M_3, s_3) \ldots$$ and the `ontological' axiom of choice, which Badiou interprets as a metaphor of ontological intervention \textit{par excellence}, ensures that one could always extend this chain through transfinite recursion beyond any given index ordinal, ultimately giving an injection from this larger ordinal into $\{(M, s)\}$, which would be a contradiction. Of course, this entails the poset  $\{(M, s)\}$ to be a set, which is true in a small category (e.g. in a Grothendieck-topos). 

Now, to establish the internal Axiom of Choice (AC), the argument again requires a proof by contradiction: assume that the maximal section $M$ is not the object $I$ itself, i.e. $M \neq I$. Because $\Sub(I)$ is Boolean, $\neg M \neq 0$. Moreover, because the subobjects of the terminal object $1$ generate the topos under the support generator (SG) condition\footnote{For an explicit definition of this condition, see my previous article, Veilahti, `Alain Badiou's Mistake', 2013.}, there must exist a non-zero morphism $t: V \to \neg M$, where $V \hookrightarrow 1$. 
The pullback of the epimorphism $X \to I$ via $V \to \neg M \hookrightarrow I$---again by an argument employing the (SG) property\footnote{See  Mac Lane \& Moerdijk, \textit{Sheaves in Geometry and Logic}, 1992, p. 277.} and a further refinement to a non-zero $W \subset V$---leads to an extension of the section $M \to X$ over the coproduct $M \cup t(W)$. This directly contradicts the initial maximality of $M$; hence, $M = I$. 

\end{proof}

When a topos satisfies both the internal AC and SG conditions, like the ones described in the previous lemma, their internal universe approximates (a bounded version of) classical set theory\footnote{To approximate the standard ZFC axiomatics, one can reconstruct the cumulative hierarchy of ordinal numbers $V_{\alpha}$ internal to the topos, as demonstrated by Fourman (see Fourman, Michael P. (1980),  `Sheaf models for set theory', \textit{Journal of Pure Applied Algebra} 19. pp. 91--101). Even under this expansion, the axiom of comprehension---the classical postulate stating that logical formulas uniquely determine subobjects---appears in a restricted configuration: for every object $b$ and formula $\phi(x)$ in which \textit{all quantifiers are restricted}, there exists a well-defined subobject $\{x \in b \mid \phi(x)\}$. Quantifiers are thus limited to bounded expressions of the type $\forall x \in b$ rather than unbounded formulations $\forall x$, aligning the framework with the `bounded Zermelo' or \textsf{BC} axiomatics.}. This approximation holds for any Cohen poset $\cond$ consisting of the finite approximations of a structural mapping $\delta \times N \to \mathbf{2}$, which, upon the successful execution of the proof, yields the required strict monomorphism $\delta \into \Power N$.

\begin{lem}
For a given Boolean topos $\Es$, to construct a strictly two-valued model in which the subobject classifier $\Omega$ possesses exactly two global elements $1 \to \Omega$, the ultra-power construction can be applied to the maximal filter $\U$ of the subobjects of the terminal object 1.
\end{lem} 

\begin{proof}
Since $\Es$ is a Boolean topos, it follows that for any subobject $U \into 1$, either the statement $U \in \U$ or its complement $\neg U \in \U$ is valid. Any arrow $1 \to \Omega$ in the filter quotient topos $\Es/\U \cong \Es_{\infty}$ corresponds to a localized arrow $1 \times U \to \Omega \times U$ over some $U \in \U$. This induces a unique morphism $f: U \to \Omega$ for a subobject $U \in \U$, where $U$ can be systematically replaced by any smaller subobject $U' \into U$ without altering the resulting global section $1 \to \Omega$ in $\Es/\U$. For a subobject $V \into U$ classified by $f$, the Boolean nature of $\Es$ ensures that either $V \in \U$ or $U \cap \neg V \in \U$. This yields two distinct factorisations:
$$\xymatrix{V \ar[r] & 1 \ar[r]^{\true} & \Omega} \quad \text{and} \quad \xymatrix{U \cap \neg V \ar[r] & 1 \ar[r]^{\false} & \Omega}$$
which demonstrates that there exist exactly two such morphisms $1 \to \Omega$ within the quotient category $\Es/\U$.
\end{proof}

As the final preliminary step before the explicit independence proof, it must be demonstrated that the representable presheaf $\mathbf{y}(p)$---which associates any condition $q \in \cond$ with the set $\Hom_{\cond}(q, p)$, defined as a singleton for $q \leq p$ and $\emp$ otherwise---constitutes a  sheaf with respect to the \textit{double negation} topology.

\begin{lem}
The representable presheaf $\mathbf{y}(p): q \to \Hom_{\cond}(q, p)$ is a sheaf.
\end{lem}

\begin{proof}
The demonstration utilizes sieves operating under the $\neg\neg$-topology, which reveal a direct structural affinity to the `dominations' $D$ utilized in \textit{Being and Event}. The critical difference is that these categorical sieves are strictly downward closed and operate as localized dominations restricted up to a specific degree, rather than acting globally over the entire poset. 

Given a $D$-indexed matching family of elements $$x_d \in \mathbf{y}(p)(d) = \Hom_{\cond}(d, p),$$ its unique amalgamation---i.e. the synthetic element whose existence is guaranteed by the sheaf condition in an analogue to Badiou's postulate of materialism---is necessarily unique because $\cond$ forms a poset. 

If the condition $q$ does not explicitly extend $p$ (which would imply the failure of the inequality $q \leq p$), there exists a coordinate $\langle a, n \rangle$ where $p(a, n) \neq q(a, n)$, or at least $q$ can be extended to a refinement $q'$ that ensures this incompatibility. This implies that $q' \leq q$, and there exists a sieve element $d \in D$ such that $d \leq q'$. However, since $p$ and $q'$ are mutually incompatible, the condition $d \leq p$ fails, directly contradicting the requirement that the family $x_d$ is matching. Thus, one must have $d \leq p$ for all $d \in D$, which proves that $\mathbf{y}(p)$ constitutes a sheaf within the $\neg\neg$-topology based on the so-called separativity of $\cond$.
\end{proof}

Now, the central theorem establishing the independence of the continuum hypothesis within the topos-theoretic framework can be formulated.

\begin{thm}
Let $\Ss$ be a two-valued Boolean topos satisfying the axiom of choice (AC) and equipped with a natural number object $N$. There exists an extended two-valued Boolean topos $\Ss'$ satisfying (AC) with a natural number object $\gamma^*(N) = \hat N$ that forces the strict inequality:
$$\hat N < \delta < \hat{\mathbf{2}}^{\hat N},$$
where $\delta = \gamma^*(\mathbf{2}^N)$ represents a large cardinality transported from the original base situation $\Ss$. Consequently, the existence of the object $\delta$ falsifies the continuum hypothesis.
\end{thm}
\begin{proof}
Because $\mathbf{2} \in \Ss$ is the subobject classifier of $\Ss$, the mapping $$f: \delta \times N \to \mathbf{2}$$ corresponds to a subobject of $\delta \times N$. Let me consider the diagonal functor $\Delta$ which embeds every object of $\Ss$ into a constant presheaf; for instance, $\Delta(N)(p) = N$ for all conditions $p \in \cond$. The objective is to construct a suitable subobject of $\delta \times N$, not directly within the base topos $\Ss$ where it is impossible due to Cantor's argument, but within the `transphenomenal' extension $\Ss \to \Sh(\cond, \neg\neg)$ expanding the exceptional algebra $\mathbf{2}$. 
The extension is nothing other than the associated sheaf functor of the canonical diagonal map, which maps each subobject $X \in \Ss$ to its corresponding constant presheaf $\Delta(X)(p) = X$ for all $p \in \cond$. Consequently, one needs to localize a specific subobject within $\Delta(\delta \times N)$ that corresponds precisely to the strict inequality one is seeking to establish.

The required subobject is first constructed at the level of presheaves as:
$$A(p) = \{(a, n) \in \delta \times N \mid p(a, n) = 0\},$$
which is the analogue of the structural subobject Badiou deploys in his demonstration. It must be shown that $A$ is closed under the $\neg\neg$-topology, that is, $\neg\neg A = A$. According to the point-wise formula \eqref{doublenegationsheaf} established above\footnote{See also Mac Lane \& Moerdijk, \textit{Sheaves in Geometry and Logic}, 1992, p. 280.}, the closure expands to:
$$\neg\neg A(p) = \{(a, n) \in \delta \times N \mid \forall q \leq p, \ \exists r \leq q: (a, n) \in A(r)\}.$$
This is the crucial step where the double negation topology on the poset becomes `transphenomenally' visible: one witnesses sequences of the form $r \leq q \leq p$ where the negation operator conjugates, operating analogously to the dense dominations $D$ that Badiou employs in \textit{Being and Event}. 

Therefore, if a pair $(a, n) \notin A(p)$, it implies that either $p(a, n) = 1$ if defined, or $p$ can be extended to a refinement forcing this value since $r \leq p$. Since every refinement $r \leq p$ preserves existing values, either $p(a, n) = 1$ or $r(a, n) = 1$, which in both cases ensures that $(a, n) \notin \neg\neg A(p)$. This demonstrates the structural inclusion $\neg\neg A \subset A$, while the reverse inclusion $\neg\neg A \supset A$ holds trivially by the axioms of closure.

Having established that $A$ is a stable $\neg\neg$-sheaf, one defines $\Omega_{\neg \neg}$ as the equalizer of the identity and double negation morphisms: $\Omega \rightrightarrows \Omega$. The characteristic map of $A$ factors uniquely through $\Delta \delta \times \Delta N \to \Omega_{\neg\neg} \to \Omega$, and the exponential transpose of this mapping yields a monomorphism of presheaves:
$$\Delta(\delta) \to \Omega_{\neg\neg}^{\Delta(N)}.$$
This is a monomorphism point-wise at the level of $\Delta(\delta)(p) \to (\Omega_{\neg\neg}^{\Delta(N)})(p)$, where the latter designates the set of natural transformations $$\mathbf{y}(p) \times \Delta(N) \to \Omega_{\neg\neg}.$$

Letting $a: \Ss^{\cond^{op}} \to \Sh(\cond, \neg\neg)$ denote the associated sheaf functor, it remains to prove that the corresponding mapping of sheaves $m: \hat \delta \to \Omega_{\neg\neg}^{\hat N}$ constitutes a monomorphism. Since the sheafification functor $a$ is left-exact, it strictly preserves monomorphisms. One must therefore verify that $a(\Omega_{\neg\neg}^{\Delta(N)}) \cong \Omega_{\neg\neg}^{a\Delta N}$. This isomorphism follows by applying Yoneda Lemma and the adjoint properties of the left-exact functor $a$, yielding\footnote{See Corollary 5 in Mac Lane \& Moerdijk, \textit{Sheaves in Geometry and Logic}, 1992, p. 282.}:
\[
\begin{aligned}
\operatorname{Hom}
(X,\Omega_{\neg\neg}^{\Delta N})
&\cong 
\operatorname{Hom}
(a\Delta N\times aX,\Omega_{\neg\neg})\
&\cong
\operatorname{Hom}
(X,\Omega_{\neg\neg}^{a\Delta N})
\end{aligned}
\]

In the resulting Boolean topos, the closed subobject classifier satisfies the isomorphism $\Omega_{\neg\neg} \cong 1 + 1 = \hat{\mathbf{2}}$, and the category possesses a natural number object $\hat N$, while $\Omega_{\neg\neg}^{\hat N}$ characterizes the internal power set $\hat{\mathbf{2}}^{\hat N}$ formed internally in the Cohen topos (after making it two-valued by the ultra-power construction). One thus arrives at the pivotal inequality $$\hat N < \widehat{\mathbf{2}^N} < \hat{\mathbf{2}}^{\hat N}.$$ Here $ \widehat{\mathbf{2}^N}$, which denotes the sheafification of the power object formed in the base topos and which Mac Lane and Moerdijk\footnote{Mac Lane \& Moerdijk, \textit{Sheaves in Geometry and Logic}, 1992, p. 283.} refer to as the `fake power set', captures the essence of the independence proof. 
 
To complete the proof, it must be demonstrated that these inequalities are strictly non-isomorphic. In any arbitrary topos $\Es$, there exist specific subobjects $\Epi(W, Z) \into Z^W$ that measure when such cardinality inequalities are strict. The structural requirement dictates that for $W < Z$, there can exist no epimorphisms---which in the set-theoretic language of consistencies translates into the non-existence of surjections $W \to Z$---for otherwise $Z$ would possess at most the same internal capacity of points as $W$. 

Without resorting to a purely diagrammatic construction of the object $\Epi(W, Z)$ analogous to the set of surjections, the definition rests on the evaluation of parameter objects of the type $\Es(E, Z^W)$ mapping arrows from $E$ to $Z^W$. In the internal Mitchell--Bénabou language---a syntax emerging relationally from the internal diagrammatic geometry of the topos (see Definition \ref{mitchell_benabou})---the object $\Epi(W, Z)$ is formalized as:
$$\Epi(W, Z) = \{f \in Z^W \mid \forall y \in Z, \ \exists x \in W: \, f(x) = y\}.$$
The demonstration concludes by showing that the triviality of $\Epi(W, Z)$ within the base model $\Ss$ necessarily forces the condition $\Epi(\hat W, \hat Z) = 0$ within the resulting Cohen topos\footnote{ Mac Lane \& Moerdijk, \textit{Sheaves in Geometry and Logic}, 1992}. The preservation of these strict cardinal inequalities follows from the Souslin property of the Cohen topos. In general, the Souslin property means that every family of pairwise disjoint subobjects is at most countable, whereas in the context of Cohen's poset this translates to any set of pairwise incompatible forcing conditions being necessarily countable\footnote{See Proposition 6 and Lemmas 7 and 8 in ibid., pp. 288--290.}. 
\end{proof}

\begin{rmk}
Unlike the classical Cohen construction reviewed by Badiou, the topos-theoretic framework bypasses the necessity of isolating a specific, `generic filter' $\female$---a despotically determined, \textit{incorporeal} point of truth whose existence remains undecidable by set-theoretic predicates alone. The category-theoretic expression of this epistemological rupture is fundamentally more abstract, or rather structural, requiring only the specification of the forcing poset $\cond$ alongside the procedural, diagrammatic filter-power construction that internalizes the `generic' procedure \textit{relationally}---the situation avoids the need to explicate a specific instance of the `generic' point $\female$. 

This non-reductive configuration is, above all, also operationally more `generic': it operates independently of point-wise ZFC axiomatics and requires no direct recourse to the transfinite ordinals of Cohen's original proof. Instead, it demonstrates how these emergent points of truth materialise a new, intensely `diagrammatic' territory. There exists a vast multiplicity of such territories, populated by contradictory points of truth. Every individual Cohen topos $\Sh( \cond, \neg\neg)$ constitutes merely a localic, and thus strictly \textit{local}, snapshot of the inherent geometry of set-theoretic expression---a `meta-ontological caesura', or rather discourse, whose very incompleteness reveals an ultimately torsional reference to subjectivity.

A specific Cohen topos, in effect, represents a particular localic configuration of Badiou's `democratic materialism', whereas alternative 
topoi---emerging from strategically distinct categories $\cond$---subjectify a multitude of diagrammatic 
effects that can never be `democratically'\footnote{In Badiou's political 
epistemology as outlined in \textit{Logics of Worlds}, `democratic materialism' 
posits a multiplicity of bodies and languages operating under local, 
consensus-driven formalism without the imposition of a singular, transcendent 
Truth. By mapping this onto topos theory, alternative 
forcing posets $\cond$ generate mutually incompatible localic universes. 
These distinct territories cannot be `democratically' unified 
or co-represented within a single, coherent global locale, which demonstrates the significance of the category-theoretic mode of expression for reframing Badiou's dialectically materialist philosophy. } 
held together within a single, isolated locale.
\end{rmk}

\section{Sheaves: Bridging Materiality and Singularity} \label{mitchellbenabou}

Badiou's philosophical project, initiated in \textit{Being and Event}, sought to mobilize mathematical language as a methodological access point to philosophize on the structural caesura as such concealed by the `ontological' prevalence of $\sets$. It was for this reason that he subsequently recognized the salience of topos theory to his account; however, he ultimately failed to execute these connections fully, restricting his philosophical updates to localic frameworks rather than adopting a genuinely topos-theoretic terminology. 

The previous sections have reconstructed these forcing arguments in a category-theoretic framework.
The importance of categorical geometry---such as the sheafification $\Ss \to \Sh(\ce, \neg\neg)$---is reflected by the fact that it actualizes \textit{novel} ways of expression, contributing to a new perspective on what it is to make a point, whether or not it acquires a distinct, ontological label (like `$\female$'). Moreover, these geometric points, as Cohen's case illustrates, can be chased from within the `ontological' topos $\Ss$, even if they were not universally valid, but conditional possibilities.

The importance of points is reflected in many domains of mathematics, such as in algebraic geometry where the Zariski-site is naturally coarse, having far too few points for topologists' taste. Capturing finer cohomological invariants precisely necessitated the introduction of Grothendieck topologies and étale morphisms to engender such new points of expression, serving as ontological anchors of otherwise invisible geometry\footnote{Étale sites are synthetically engineered to be significantly more `tensed' and structurally complex than the traditional Zariski-site of algebraic varieties, providing a richer environment in which previously invisible geometric and cohomological structures can be expressed.}. Indeed, a `weakly' or `dialectically materialist' topos, understood here as a Grothendieck topos of set-valued sheaves, is an interesting boundary case, where internal geometries communicate with external representations. In fact, even if a Grothendieck topos does not have enough points to begin with---especially in the `atonic' case---Johnstone proves that there exists a unique smallest topology $j$ such that the subtopos $\Sh_j(\Es)$ has 
'enough points'\footnote{See Johnstone, \textit{Topos Theory}, 1977, p. 242.}. This fact might be viewed as supporting Badiou's external, `ontological' approach to topology---yet it obviously does not apply when one shifts to non-material topoi, which cannot be ontologically represented as sheaves of sets.

Regardless of whether a topos is material (over sets), the remaining question relates to the philosophical meaning of these forcing procedures, especially when interpreted internally to a topos, rather than as externally imposed, model-theoretic conditions. At the same time, I will link this \textit{problem of semantics} and internal forcing conditions to sheaf theory. This serves two purposes. First, it helps to make the `torsion' inherent to the internal concept of a topos materially meaningful, in a way appealing to what \BE\ refers to as `deductive fidelity'. Second, the internal reframing of semantics further demonstrates the limitations inherent to the localic worlds in \LW, and why a true connection between the  two works is conceivable only by abandoning the strong, atomistic condition of materiality.

So what do points, which can be defined either geometrically or as binary conditions, fundamentally stand for? Classical analytic philosophy maintains a strictly \textit{propositional} relation to points of truth. To make the geometric definition of points meaningful, it is necessary to bridge these relationally emergent, sheafified constructs with an \textit{internal} representation of deductive propositions. 

I already observed during the previous proof the need for a systematic \textit{translation} between the `internal' 
or structural phenomena of a topos and the language executed on the 
operational level of proofs---a translation that subsequently gives propositional logic 
and language a dual status: two modes of articulation related to the distinction between designation and expression.

\BE\ encounters this dilemma in an inexplicated, sometimes confusing way.  \LW, while perpetuating this ambiguity, seeks to demonstrate that the local languages of a topos are not fundamentally different from the `ontological' reality of sets---a claim that he actually \textit{forces} to be true by imposing the localic condition on a topos (i.e. the strong, atomistic postulate). Badiou's expression of the forcing conditions, which undoubtedly forms a core philosophical contribution of \BE, thus becomes arrested by an indeed `unconscious' ontologist behind the analytically preoccupied phenomenology in \LW. He not only seeks to enforce specific conditions of expression, but more urgently mistakes the possibility of their (set-theoretic) conditioning as being their sole condition of possibility. As a consequence, Badiou refutes any genuine event or individuation concerning topos-theoretic descriptions, simultaneously delimiting his own structural toolkit. The localic phenomenology of \LW\ thus regulates the concept of regulation itself in terms internal to a topos with its strictly atomistic presentation of truth. 

Compatible with this philosophical landscape, the Cohen topos consisting of sheaves over the set of `conditions', while `generalizing' Badiou's\footnote{See Badiou, \BE, 2006, pp. 416--421.}
 model-theoretic extension $S \into S(\female)$, marks at best a `fact' but not a true singularity, not yet fully altering the body of mathematical expression. While it structurally dictates that in a binary sense undecidable statements retain a decisive valuation within the modified situation, it is still preoccupied by the specific, `ontological' forcing poset $\cond$. Instead, for an honestly new body of expression outside Badiou's atomistic agenda over mathematical formalism, a method must be established to translate these statements into categorical terms that can receive a stable, diagrammatic truth-value irreducible to the forcing of binary points. 
 In other words, for a topos to truly carry the burden of meta-ontological caesura capable of actualizing torsion inaccessible to the split language of sets, and thus under the interdict of Badiou's fundamental law, one needs a body where this language itself is merely a local, contextual effect, rather than a causal basis of its construction. For it to truly make a singular expression surface, it is mandatory to reverse the translation vector between category theory and set theory, conceiving the `ontological' architecture itself in strictly \textit{internalized} terms. 

Crucially, Cohen's original proof is conducted within classical set-theoretic language. But this language itself can be formalized diagrammatically as a language operating upon the \textit{internal} Heyting algebra $\Omega$---a structure fundamentally distinct from the external `transcendental indexing' employed by Badiou\footnote{Most strikingly, when Badiou claims that $T$ retains `degrees of 
intensity', these degrees are precisely \textit{extensive} and not intensive, given how it emerges as a materialist projection of a Grothendieck topos (i.e. a topos that satisfies the weak postulate of materialism). 
Badiou obviously targets this philosophical critique against Gilles Deleuze's discussion 
of the `intensive' and `extensive' dimensions of individuation. See Deleuze, Gilles (1994), \textit{Difference and Repetition}. Transl. Paul Patton, New York: Columbia University Press. [Originally published in 1968.] Discrediting Badiou's philosophical 
argument based on his atomistic vision of topos theory---which artificially 
bridges the internal expression and the language designated to mirror operational dogmatism---it becomes obvious that the internal `grading' of the subobject classifier 
algebra $\Omega$ is infinitely more `intensive' than the external, localic grading $T$.}.
 Therefore, for Cohen's proof to remain `internally' consistent with the `phenomenon' of points \textit{within} a topos, it must be effectuated internal to the topos, instead of serving as its causal condition. This path is realized by utilizing the internal diagrammatic syntax of the \textit{Mitchell--Bénabou} language: 

\begin{dfn}[Mitchell--Bénabou Language]
The \textit{Mitchell--Bénabou} language $L_{\Es}$ of a topos $\Es$ consists of the following structural components:\footnote{See Johnstone, \textit{Topos Theory}, 1977, p. 153. See also Mac Lane \& Moerdijk, \textit{Sheaves in Geometry and Logic}, 1992, pp. 298--300.}
\begin{enumerate}
\label{mitchell_benabou}
\item \textbf{Types:} The types of the language are the objects $X, Y, Z, \ldots$ internal to the topos $\Es$.
\item \textbf{Terms:} The terms are recursively generated and belong to specific types:
\begin{itemize}
\item There exists an infinite stock of \textit{variables} $\mathbf{x}, \mathbf{x'}, \mathbf{x''}, \ldots$ of each type $X$.
\item For each arrow $f: X \to Y$ in $\Es$ and any term $\tau$ of type $X$, there is a corresponding term $f(\tau)$ of type $Y$.
\item For any two terms $\sigma$ of type $X$ and $\tau$ of type $Y$, there is a product term $\langle\sigma, \tau\rangle$ of type $X \times Y$.
\end{itemize}
\item \textbf{Formulas:} The formulas in the language $L_{\Es}$ are defined as terms whose type is the subobject classifier $\Omega$---i.e. internal expressions of transcendental `truth-values':
\begin{itemize}
\item For each type $X$, there is a binary equality predicate $=_X$ of signature $(X, X)$, where the formula $(\sigma =_X \tau)$ is defined via the characteristic morphism $\delta_X: X \times X \to \Omega$ classifying the diagonal embedding\footnote{The equality formula $(\sigma =_X \tau)$ corresponds structurally to the composition $\delta_X \circ \langle\sigma, \tau\rangle$.}.
\item For each type $X$, there is a binary membership predicate $\in_X$ of signature $(X, \Omega^X)$, defined internally by the canonical evaluation morphism $X \times \Omega^X \to \Omega$.
\item The logical connectives $\wedge, \vee, \Rightarrow$, and $\neg$ are induced directly by the internal Heyting algebra operations acting upon $\Omega$.
\end{itemize}
\item \textbf{Quantifiers and Descriptions:} If $\phi$ is a formula and $\mathbf{x}$ is a variable of type $X$, then:
\begin{itemize}
\item $\forall \mathbf{x} \phi$ and $\exists \mathbf{x} \phi$ are formulas of type $\Omega$.
\item $\{\mathbf{x} \mid \phi\}$ is a term of type $\Omega^X$, mapping the internal analogue to the axiom of comprehension.
\item $\iota \mathbf{x} \phi$ is a description term of type $X$, denoting the unique element satisfying $\phi$ if such an element exists. This is a crucial instrument for the internal expression of both the axioms of comprehension and separation, and it remains equally necessary for the internal expression of forcing arguments, such as Cohen's construction.
\end{itemize}
\end{enumerate}

As in Badiou's localic world enforcing the atomistic projection of $\Omega$, also in a general elementary topos---or what I have come to label a singular topos, as opposed to the `higher' compositions that prevent ontology or the `fundamental law' appearing even internally---the law of \textit{reductio ab absurdum} does 
not hold. Consequently, mathematical interventions must restrict 
themselves to constructive, intuitionist arguments, since the double negation 
topology $\neg\neg$ does not necessarily coincide with the identity\footnote{Mac Lane \& Moerdijk, \textit{Sheaves in Geometry and Logic}, 1992, pp. 298, 302.}. While the logic Badiou constructs upon the external transcendental frame is similarly intuitionist, it remains external to the topos; its points are invariably subjected to a split, classical choice of truth: either/or. Propositions in Mitchell--Bénabou language instead relate strictly to the internal body of transcendental representation $\Omega$, being intensive to the extent that it retains no simple, element-wise representation in classical, set-theoretic terms.
\end{dfn}

\begin{rmk}
Badiou's engagement with locales could obviously be read as an internal 
manifestation of his supposedly non-localic topos, enacted specifically on the internal transcendental locale $\Omega$. That interpretation 
would, however, introduce profound philosophical, rather than purely mathematical, problems. Even the `constitutive' materialism 
postulated by his atomistic, `real synthesis' would reduce to a mere structural effect, stripped of any causal necessity 
regarding the topos as an ontological, privileged construct. In particular, it would 
effectively invalidate any phenomenological argument concerning the limits 
of the structured appearance of being. Badiou's 
propositionally performed claim about `materiality' would thus constitute nothing 
more than a linguistic expression, positioning himself in the relational world that he specifically criticizes. Ultimately, he fails to `ontologically' constrain either the conditions 
of weak materiality or the claimed atomism of the topos---features that are certainly `internal' realities of a singular (elementary) topos, which in my view is the natural home of `dialectical materialism'.
\end{rmk}

\begin{rmk}
For an internal language to `make sense'---in a manner analogous to the interpretation function that provides the material coordinate of subjectivity in \textit{Being and Event}---it requires the specification of an `interpretation' functor operating strictly in terms of the morphisms of $\Es$\footnote{A type $X$ is interpreted as the identity morphism $\id_X: X \to X$. A term $f(\tau)$, where $\tau$ is of type $X$ with a stock of free variables $(\mathbf{u}_i)_{i=1}^n$, is interpreted as the composite morphism $f \circ \lvert\tau\rangle: \prod U_i \to X \to Y$. For instance, the comprehension term $\{\mathbf{x} \mid \phi\}$ evaluates to a morphism $\lvert\{\mathbf{x} \mid \phi \}\rvert: \prod U_i \to \Omega^X$, which is precisely the exponential transpose of the structural interpretation morphism $\lvert\phi\rvert: X \times \prod U_i \to \Omega$. In this sense, the Mitchell--Bénabou language systematically translates terms into subobjects (morphisms targeting $\Omega$), where stable `truth'-values---that is, closed `statements'---reduce to subobjects of the terminal object, which is to say, global sections of the form $1 \to \Omega$. Concurrently, a formula containing free variables, for example of type $X$, evaluates directly to a subobject of $X$, classified by an arrow $X \to \Omega$.}.

To illustrate how the universal quantifier $\forall$ operates internally, consider a formula $\phi$ of signature $(X, Y)$ containing the free variables $\mathbf{x}$ and $\mathbf{y}$. The quantified formula $\forall\mathbf{x} \phi(\mathbf{x}, \mathbf{y})$ is no longer free with respect to the variable $\mathbf{x}$, and consequently, its semantic interpretation must be mediated by an arrow of the form $Y \to \Omega$. Given the unique terminal morphism $p: X \to 1$, there exists an induced pull-back frame morphism $\Power p: \Omega^1 \to \Omega^X$, which possesses the internally defined adjoint functors $$\forall_p  \textrm{ and } \exists_p: \quad \Omega^X \to \Omega^1 \cong \Omega.$$ Since the formula $\phi(\mathbf{x}, \mathbf{y})$ expresses a generalized term $Y \to \Omega^X$, its composition with the internal right adjoint $\forall_p$ yields the desired structural morphism $Y \to \Omega$\footnote{Ibid., p. 300.}.

Recall that truth is conceptualized as a specific global subobject of the terminal object, $1\lhook\joinrel\xrightarrow{\text{true}} \Omega$. Therefore, a statement expressing a truth within the topos is structurally defined as one whose semantic interpretation factorization forces the emergence of this specific morphism. If a formula $\phi$ with a framework of free variables $(\mathbf{u}_i)$ induces a composite factorization of the form:
$$\text{true} \circ p: X \times \prod U_i \to 1 \xrightarrow{\;\true\;} \Omega,$$
the formula $\phi$ is said to be \textit{universally valid} within the topos, a structural state conventionally denoted by the forcing notation:
$$\Es \models \phi.$$  
\end{rmk}

\begin{rmk}
It is equally worth noting, with respect to basic logical calculus, that if a topos possesses a natural number object $N$, a vast array of counting operations and `ontological' rules governing them remain valid internally to the topos. Provided that proofs are strictly constructive, they are expressed by a topos, as are all finitary algebraic theories. 

However, one is not required to restrict the investigation to mere algebraic syntax. There is an internal theory of functions analogous to the Mitchell--Bénabou language, including (1) a stock of internal types $X, Y, Z, \ldots$; (2) a family of \textit{function symbols} $\alpha, \beta, \ldots$ governed by distinct signatures; (3) a collection of \textit{primitive relation symbols} $r, s, \ldots$; (4) a recursive generation of \textit{terms} encompassing variables $\mathbf{x}, \mathbf{x'}, \ldots$ alongside the evaluation of function symbols over existing terms $\alpha(\bar{\tau})$; (5) a framework of \textit{formulas} comprising the atomistic constants $\true$ and $\false$, internal equality $=$, relation evaluations $r(\bar{\tau})$, geometric conjunctions and disjunctions, and the existential formula $\exists \mathbf{x} \phi$; and, finally, an internal representation of categories\footnote{See Johnstone, \textit{Topos Theory}, 1977, pp. 199--200.}.
\end{rmk}

\section{Kripke--Joyal Semantics and Badiou's Connection to Fidelity}

Topos theory replaces not just the traditional, binary composition of points---it alters the topology of truth itself. 
Of course, such a claim indicates its own hierarchies---it is about making a point about topos theory, after all. Yet by viewing this debate itself as a topological effect, instead of as a manifestation of ontological force, the roles of contents and expression are swapped, in a way that resembles Louis Hjelmslev's semiotic concept of double articulation\footnote{See Deleuze, Gilles \& Félix Guattari (1988), \textit{Thousand Plateaus---Capitalism and Schizophrenia}. Transl. Brian Massumi, Minneapolis, University of Minnesota Press. [Originally published in 1980.]}. This framework is 
relativistic in the sense that it denies the primacy of either modality of signification, undoing any transcendental hierarchy between them. This observation is particularly fruitful in the context of mathematics, as Hjelmslev originally further separated the form and substance of both expression and contents---yet in mathematics, at least for the formalist school, the form itself is substance. 

Topos theory further bridges this concept of double articulation by systematically `internalizing' set-theoretic expression as the diagrammatic `content' of a specific topos, thus categorically diminishing any distance between the different sides of signification. Yet it does not mean that a topos is not a meaningful concept. Quite the contrary, if a hierarchy itself---like the one endowed by points---is an effect of a topos, it becomes the very condition of meaning and expression. For meaning to acquire its formative power, \textit{semantics} is the mediator between purely relational existence and the local, contextual effects of operational force. Semantics, which is a precedent of 
narrative semiotics, thus endorses the capacity of the internal concept of order $\Omega$ to mobilize itself in a transcendental (or `dialectical') way amidst the apparently plural logics of worlds. Yet, unlike in an `atomistic' topos, this power emerges from within, rather than as a necessary or ontologically prescribed structure imposed `from above'. 

That localizing effect of force internal to a topos has already swallowed Badiou\footnote{Badiou, \textit{Being and Event}, 2006, p. 358.} as he turns to \textit{Tractatus Logico-Philosophicus}, observing that while Wittgenstein `still only admits the existence of a single world', he already `saw perfectly well that the general form of a world, or that which is intelligible in the worldliness of a world, is purely logical'. Even while a topos is a different animal, it has the power to lure oneself into submitting to the local force of language. Therefore, a topos does not merely topologize spatial and temporal coordinates, but the very \textit{semantics of truth itself}. 

Precisely if one follows Wittgenstein\footnote{Wittgenstein, Ludwig (1981), \textit{Tractatus Logico Philosophicus}. Transl. C. K. Odgen. Oxon and New York: Routledge. [Orig. published 1922.] Prop. 4.45.} in asserting that `to give the essence of a proposition means to give the essence of all description, and thus the essence of the world', the diagrammatic `essence' of semantic expression has already displaced the world constituted by the ontological, split realm implicated by the question of consistency (which privileges contents). Therefore, Wittgenstein's fidelity in the essence of language requires the structure of truth itself to appear \textit{topo}-logically. It aims to respond to this precise dilemma: how can the world be simultaneously plural and infinite, and yet produce the effect of finitude that grounds its localic and logical perception. 

To capture this dilemma, the crux is to recognize that the subject of forcing itself---the very paragon of semantic expression---is itself a spatio-temporal, topo-logical process. Cohen's original proof, as Badiou traces it, was grounded upon the mechanism of \textit{forcing}; the genius of topos theory was then to \textit{contextualize} this very mechanism and to force the apparent incompleteness of language to reveal itself not as a constitutive reality of being, but instead as an effect of that transcendental, localizing force itself when thought of as if imposed from outside. A topos provides a direct interpretation of the dilemma by mobilizing a contextual understanding of deductive and ultimately `ontological' subjectivity. To express this formally, let me now introduce the passage from Kripke's original work on the semantics of forcing to its topos-theoretic reinterpretation as a contextual, localized effect of force. 

\begin{dfn}[Forcing Semantics]
The forcing relation $U \models \phi(\alpha)$ is uttered as `$U$ forces $\phi(\alpha)$' within a topos $\Es$, where $\alpha: U \to X$ constitutes a generalized element of an object $X$, governed by the semantic rules regulating when $\alpha$ structurally `belongs' to a subobject designated by the internal comprehension notion $\{x \mid \phi(x)\}$\footnote{In diagrammatic terms, this means that the generalized element $\alpha: U \to X$ factors uniquely through the monomorphism $\{x \mid \phi(x)\} \hookrightarrow X$, conveying that the image subobject satisfies $\Im(\alpha) \leq \{x \mid \phi(x)\}$ within the Heyting algebra of subobjects of $X$. This relation is manifested in the following commutative diagram:
$$\xymatrix{
& \{x \mid \phi(x)\} \ar@{^(->}[d] \ar[r] & 1 \ar[d]^{\true} \\
U \ar@{-->}[ur] \ar[r]_{\alpha} & X \ar[r]_{\phi} & \Omega
}$$
If the formula $\phi$ contains no free variables, the forcing statement $U \models \phi$ demonstrates that the unique terminal morphism $U \to 1$ factors through the global subobject $\{\;\cdot \mid \phi\} \hookrightarrow 1$ classified by the closed section $\phi: 1 \to \Omega$.}. Specifically, $1 \models \phi$ if and only if the semantic relationship reduces to the global section $1\lhook\joinrel\xrightarrow{\text{true}} \Omega$, carrying the interpretation that the sentence is not merely locally forced, but `globally' effective. 
The fact that this is not always the case evidences how a topos localizes not just what mathematical entities, including the arrow `true', truly \textit{are}, but that it topologizes propositions themselves, that is, the very fidelity of the subject seeking to propose. \end{dfn}

\begin{prop}[Monotonicity and Soundness]
Kripke--Joyal semantics, which topologizes the concept of forcing, is topos-theoretically well-defined in the sense that it is strictly compatible with the composition of arrows as well as the local characteristics of a topos\footnote{If $U \models \phi(\alpha)$, then $U' \models \phi(\alpha \circ f)$ for any morphism $f: U' \to U$. Conversely, if $f: U' \to U$ is an epimorphism, the implication operates in the reverse direction: $U' \models \phi(\alpha \circ f) \Longrightarrow U \models \phi(\alpha)$.}. These properties are critical for `materializing' or objectifying what Badiou regards as subjects' `deductive fidelity'; in other words, it formally bridges the theory of subjectivity in \BE\ to the \textit{sheafification process} essential to Badiou's postulate of `dialectical materialism' in \LW.
\end{prop}

\begin{thm}[Kripke--Joyal Semantics]
The internal forcing relation satisfies the following recursive conditions, conforming structurally to the naïve semantics of set theory:
\begin{enumerate}
\item $U \models \phi(\alpha) \wedge \psi(\alpha)$ if and only if $U \models \phi(\alpha)$ and $U \models \psi(\alpha)$.
\item $U \models \phi(\alpha) \vee \psi(\alpha)$ if and only if there exist morphisms $p: V \to U$ and $q: W \to U$ such that the coproduct copairing $p + q: V + W \to U$ is an epimorphism, while concurrently $V \models \phi(\alpha \circ p)$ and $W \models \psi(\alpha \circ q)$.
\item $U \models \phi(\alpha) \Rightarrow \psi(\alpha)$ if and only if for any morphism $p: V \to U$ such that $V \models \phi(\alpha \circ p)$, it follows that $V \models \psi(\alpha \circ p)$.
\item $U \models \neg \phi(\alpha)$ if and only if whenever a morphism $p: V \to U$ satisfies $V \models \phi(\alpha \circ p)$, it follows that $V \cong 0$.
\item $U \models \exists y \phi(\alpha, y)$, where $y$ is a variable of type $Y$, if and only if there exists an epimorphism $p: V \to U$ and a generalized element $\beta: V \to Y$ such that $V \models \phi(\alpha \circ p, \beta)$.
\item $U \models \forall y \phi(\alpha, y)$, where $y$ is a variable of type $Y$, if and only if for every object $V$, every morphism $p: V \to U$, and every generalized element \mbox{$\beta: V \to Y$}, it follows that $V \models \phi(\alpha \circ p, \beta)$\footnote{The universal quantifier condition is equivalent to the statement that the product domain satisfies $U \times Y \models \phi(\alpha \circ \pi_1, \pi_2)$. However, for an arbitrary object $C$ within the underlying site of a Grothendieck topos, there exists no such simple set-theoretic interpretation for the universal quantifier, demonstrating the necessity of the internal localic mapping.}.
\end{enumerate}
\end{thm}

\begin{proof}
For illustration, let me examine the proof of the first condition, leaving the less straightforward cases as verification exercises for the reader\footnote{See Mac Lane \& Moerdijk, \textit{Sheaves in Geometry and Logic}, 1992, pp. 304--305.}. Let the morphisms $\phi(x)$ and $\psi(x): X \to \Omega$ classify their corresponding subobjects $\{x \mid \phi(x)\} \hookrightarrow X$ and $\{x \mid \psi(x)\} \hookrightarrow X$. By the definition of the logical conjunction operation $\wedge: \Omega \times \Omega \to \Omega$, one obtains a canonical pull-back diagram:
$$\xymatrix{ 
\{x \mid \phi(x) \wedge \psi(x)\} \ar@{^(->}[r] \ar@{^(->}[d] & \{x \mid \phi(x)\} \ar@{^(->}[d] \\
\{x \mid \psi(x)\} \ar@{^(->}[r] & X,
}$$
which renders the assertion immediate. The internal forcing of an equality, written as $U \models \sigma(\alpha) = \tau(\alpha)$, reduces directly to establishing the identity $\sigma' \circ \alpha = \tau' \circ \alpha$, where $\alpha: U \to X$ is a generalized element and $\sigma', \tau'$ constitute the structural interpretations of the terms $\sigma$ and $\tau$.
\end{proof}

\begin{rmk}
The most critical primitive predicate of the Mitchell--Bénabou language is the membership relation $\in_X$, defined internally via the canonical evaluation morphism $X \times \Omega^X \to \Omega$. Indeed, this predicate $\in_X$ is itself structurally isomorphic to a \textit{subobject} of the product $X \times \Power X$---a relation signifying which elements belong to which specific subobject---classified precisely by the monomorphism $\in_X \hookrightarrow X \times \Omega^X$\footnote{If $\sigma$ is a term of type $Y$ and $\tau$ is a term of type $\Omega^Y$ containing a variable of type $X$, then for any generalized element $\alpha: U \to X$, the forcing of the membership relation $U \models \sigma(\alpha) \in_Y \tau(\alpha)$ is equivalent to stating that the pairing morphism $\langle \sigma'(\alpha), \tau'(\alpha)\rangle: U \to Y \times \Omega^Y$ factors through the subobject embedding $\in_Y \hookrightarrow Y \times \Omega^Y$.}.

With respect to the `internal' production of Badiou's foundational event-notation, `$e \in e$', one must interrogate how the Mitchell--Bénabou language processes such a self-referential syntax. Crucially, such an expression is syntactically unformulated given the rigid signature $(X, \Omega^X)$ of the internal membership predicate $\in_X$. The language cannot relate a term to itself under this predicate because the types of the two sides of the operator disagree (i.e. $X \neq \Omega^X$)---hence the prohibition of the event reveals a structural link to the internal expression of ontological `excess' (i.e. a quality of the power-operator) irreducible only to the axiom of foundation\footnote{See ft. \ref{axiomoffoundation}, p. \pageref{axiomoffoundation}.}.

\end{rmk}

\begin{rmk}
As a prominent example of how predicate logic becomes topos-theoretically embodied, this geometric language is fully capable of expressing the internal version of the \textit{axiom of choice} (AC), which remains constitutive to Badiou's\footnote{Badiou, \textit{Being and Event}, 2006, p. 228.} association of subjectivity with logical interventions. It states that the exponential transpose  $\Omega^X \to \Omega^Y$ of any epimorphism $f: X \to Y$ is itself an epimorphism. Formally, this reduces to the internal assertion:
$$\textrm{(AC)} \qquad \forall f \in Y^X [(\forall y \exists x \, (f(x) = y)) \Rightarrow \exists s \in X^Y \forall y \in Y \, (f \circ s(y) = y)], $$
which establishes that in terms of generalized elements, every internal epimorphism `splits'---another indication of the idiosyncrasy of his intervention-concept that folds around the binary concept of choice.
\end{rmk}

\begin{prop}[Localic Forcing]\label{localeproposition}
In the specific case of a locale---that is, a world fulfilling Badiou's strong, atomistic postulate such as the Cohen topos or more generally a topos of presheaves $\sets^{T^{\op}}$---the internal forcing semantic relation yields the conventional rules of forcing regulated via `degrees of appearance $p \in T$:\footnote{See Mac Lane \& Moerdijk, \textit{Sheaves in Geometry and Logic}, 1992, p. 318.}
\begin{enumerate}
\item $p \models \phi(\alpha) \wedge \psi(\alpha)$ if and only if $p \models \phi(\alpha)$ and $p \models \psi(\alpha)$.
\item $p \models \phi(\alpha) \vee \psi(\alpha)$ if and only if for any condition $q \leq p$ there exists a further refinement $r \leq q$ such that $r \models \phi(\alpha \circ r)$ or $r \models \psi(\alpha \circ r)$. 
\item $p \models \phi(\alpha) \Rightarrow \psi(\alpha)$ if and only if for all conditions $q \leq p$ for which $q \models \phi(\alpha \circ q)$, it follows that $q \models \psi(\alpha \circ q)$.
\item $p \models \neg \phi(\alpha)$ if and only if there exists no refinement $q \leq p$ such that $q \models \phi(\alpha \circ q)$.
\item $p \models \exists y \phi(\alpha, y)$ if and only if for any $q \leq p$ there exists a refinement $r \leq q$ and a local section $\beta \in Y(r)$ such that $r \models \phi(\alpha \circ r, \beta)$.
\item $p \models \forall y \phi(\alpha, y)$ if and only if for any condition $q \leq p$ and any local section $\beta \in Y(q)$, it follows that $q \models \phi(\alpha \circ q, \beta)$.
\end{enumerate}
\end{prop}

This atomistic, still `democratically materialist', semiotics remains fundamentally simple: it demonstrates that localic forcing theory is entirely compatible with the intuitive operations of classical logic. In other words, in Badiou's local, quasi-Platonic worlds of `calculated phenomenology', the concept of semantics remains bounded in that no real torsion is expressed: it remains a virtual idea that is not materially actualized except under materially weak, non-atomistic topoi of general sheaf theory. 

That structural annihilation of semantic torsion is not preserved in general, however. Typically, semantics 
cannot be materialized within an atomistically split, set-theoretic sense 
at all. Instead, sheafification is the precise key to unlocking how semantics itself can materialize torsion in a manner inconceivable to Badiou, that is, to assign the unassignable. As mentioned earlier, a sheaf-object replaces classical `sets of functions' within the categorical 
language of mathematics (see Remark~\ref{functions}, p.~\pageref{functions}). 

The same now applies to semantics whose `points' no longer need to depend on element-wise truth-values. Sheafified semantics thus mediates the passage between a relational truth-concept and the incorporeal, atomistic realm of materiality pertinent to Badiou's point-wise understanding of being. A dialectically materialist topos---whose grading is not necessarily localic, implying that it may sustain the idea of \textit{torsion}\footnote{Weak materialism represents the `materialist' half of Badiou's postulation of materialism. Badiou's original postulate consists, through its formal mathematical analysis, of two different requirements: that the topos is a Grothendieck topos possessing a geometric morphism $\gamma: \Es \to \sets$, and that the topos supports a set of generators. Only in the combined case is the morphism $\gamma: \Es \to \sets$ strictly bounded and logical, thereby reducing to a locale graded over the external complete Heyting algebra $\gamma_*(\Omega)$. See Veilahti, `Alain Badiou's Mistake', 2013.}---poses the question of how even an `unconscious ontologist', like Badiou's set-theorists, can make a point about torsion; even without the least claim to understanding its evental meaning.

To explore this scenario, there is a need to establish a method relating semantics contextualized upon generalized elements $U \to X$ with a \textit{sheafified semantics} similarly contextualized by an index category $\ce$ that is \textit{not necessarily a locale}. 
To this end, let me consider a Grothendieck topos over a site category $\ce$ 
equipped with a Grothendieck topology $j$. There exist canonical embedding 
and localization functors:
$$\xymatrix{
\ce \ar[r]^{\mathbf{y}} & \sets^{\ce^{\op}} \ar@{->}@<2pt>[r]^{\a} & \Sh(\ce, j) \ar@{->}@<2pt>[l]^{\iota} & (= \Es),
}$$
where $\mathbf{y}$ denotes the Yoneda embedding which maps an object of $\ce$ to its representable presheaf, and $\a$ represents the associated sheaf functor taking each presheaf to its sheafification. For any object $C \in \ce$, it is possible to evaluate the sheafification $\a\mathbf{y}C$ and consider the morphisms $\alpha: \a\mathbf{y}C \to X$. By the adjoint property $\a \dashv \iota$ and the application of the Yoneda lemma, such a morphism corresponds uniquely to a specific element within the set $X(C) \cong \Hom_{\Es}(\a \mathbf{y} C, X)$. 

Therefore, not only generalized elements $U \to X$, but also the objects of the underlying indexing site $\ce$, can be treated as explicit conditions of forcing, written as:
$$C \models \phi(\alpha).$$ 
This is structurally equivalent to the statement that $\alpha \in \{x \mid \phi(x)\}(C)$, or alternatively, that the morphism factors through the structural embedding $\mathbf{y}C \to \{x \mid \phi(x)\} \hookrightarrow X$. For example, if one assigns $\ce = \Cond$, it can be examined for every condition $C \in \cond$ the explicit local section set $X(C)$, where the formula $\phi(x)$ isolates the sub-sheaf relation:
$$\{x \mid \phi(x)\}(C) \subset X(C).$$

\begin{dfn}[Sheafified Semantics]
The sheafified semantics $C \models$, defined canonically over the site category, 
satisfies the same monotonicity and local characteristics as the generalized 
forcing relation $U \models$; for instance, $C \models \phi(\alpha)$ implicates $C_i \models \phi(\alpha \circ f_i)$ for any site morphism $f_i: C_i \to C$. 
The algebraic properties governing $C \models$ correspond structurally to 
those of $U \models$, with the explicit exception of the quantifier modalities:
\begin{enumerate}
\item $C \models \exists y \phi(\alpha, y)$ if and only if the object $C$ has a suitable $j$-covering family $\{f_i: C_i \to C\}_{i \in I}$ and a corresponding family of local sections $\beta_i \in Y(C_i)$ such that $C_i \models \phi(\alpha \circ f_i, \beta_i)$ holds for each index $i \in I$. 
\item $C \models \forall y \phi(\alpha, y)$ if and only if for every morphism $f: D \to C$ and every generalized local section $\beta \in Y(D)$, it follows that $D \models \phi(\alpha \circ f, \beta)$. 
\end{enumerate}
\end{dfn}

\begin{rmk}
With such a sheafified, weakly materialized semantics---which can be held as the paragon of dialectical materialism---the verification of the independence of the continuum hypothesis reduces to the exact procedure originally executed by Cohen over $\sets$, yet crucially circumventing any explicit, external reference to Badiou's non-constructive `generic multiple' $\female$. The `unassignable' or the `inexistent' spreads across an entire topological landscape, without the need for incorporating it through explicit, elementwise names like `$\female$'.  

Furthermore, by interpreting the recursive properties of forcing semantics 
directly upon this category of sheaves, it becomes apparent that the closed 
formulas corresponding to global sections $1 \to \Omega$ function precisely 
as the `points of truth' of the topos, corresponding to subobjects of the classifier $\Omega$, but again without an `ontological', elementwise description. It is only when the `transcendental' topological indexing $\Omega$ is \textit{forced} into the frame of a poset that, given a formula $\phi(\alpha)$, 
the collection of conditions $\{p \in T \mid p \models \phi(\alpha)\}$ 
amounts to a sub-Heyting algebra of the forcing poset $T$, effectively inducing a canonical homomorphism into the binary lattice $\Two$ and directly adopting the elementwise description of a point depicted in \LW.

As I asserted above, within the `ontological' ground category $\Ss$, these statements remain \textit{ipso facto} `undecidable' and thus fail to constitute points of truth; it is precisely 
the sheafification over the $\Cond$-poset that effectuates their `material' possibility---an act or process \textit{external} to the situation itself that procedurally adopts the subjectifying role implicated by the `fundamental law' in \BE. Therefore, this moment central to `dialectical materialism' specifically reveals  
another layer of the concept of forcing itself, which is structural in nature rather than being bound to some specific incorporeal name $\female$. Instead, it marks the very motif of localization itself, which both illustrates how the `dialectical' method itself constitutes a subjective operation of force, that in turn makes it possible to ask the question: whether or not a topos, as a whole, is atomistic or localic (i.e. whether the transcendental object $\Omega$, or the conditioning category $\ce$ over which sheaves are `materialized', is an external poset). 

Therefore, what Badiou presents as an alternative to the local condition of mathematical language is, in fact, a localizing force in another expressive register: the register of relational descriptions. While topos theory is fully capable of formalizing classical set theory in  
the internal syntax of the Mitchell--Bénabou language---and thus without the need to force an elementwise, `ontological' translation of points themselves---it at the same time effectuates another question about the adequacy of materiality itself and whether the `internal atom' $\Omega$ can again be reterritorialized by an external indexing category $\ce$, at the same time subjecting the world to another question: `dialectical' or `democratic'. But this question is by no means immanent to the topos itself, outside its own, internal self-description. Yet Badiou, revealing his materialist aspirations towards `Platonic' ontology outside of what he can envision \textit{only} the singular event or the inexistent, can only understand this world from within. Taken up by its internal semantics, he cannot perceive an alternative to its `fundamental law' except through its occultation, thus unable to escape this fundamental question between atomism and dialectics by stepping outside the singular, `dialectical' topos entirely.

\end{rmk}

\section{What's That `Another Event'? Temporalizing Badiou's Gap }

Topos theory  not only reveals the dual characteristics of logical language, but it also poses direct challenges to the tradition that posits set theory as a mathematical meta-ontology and the locus of expression \textit{par excellence}. Undeniably, proofs still demand a point. Yet this is not equal to saying that structural compositions, like a topos, should be made causally liable to the existence of a proof. Better yet, they can internally express a process of deductive inquiry, regardless of whether one assumes their force in the domain of ontic necessity. 
Thus it remains pivotal to examine how the structural framework within which Badiou theorizes the \textit{enquiry} of the local configuration of the event intersects with topos theory. 

Badiou himself denotes this enquiry as a specific, local and finite process through which the consequences of the `event'---i.e. of novel mathematical expressions as in my case---are systematically \textit{verified}. Crucially, this truth-process is temporal: it does not emerge with the same computational instantaneity that characterizes the traditional, idealist composition of bi-valued set-theory. 
Therefore, in spite of Badiou's bifurcated, split understanding of mathematical regulation, he poses a concept of duration---a `\textit{durée}'---though only presenting it in an abrupt form: a `gap' that separates the `evental multiple-being' $e$ from the corporeal connection of fidelity, that is, from the domain of propositions. 

In terms of formal notation, Badiou marks this `historical' enquiry with the predicate $\square$, which is existentially sluggish in the sense that it can be determined only across a finite succession of cases, directly opposing the `statist' thesis of the classical membership predicate $\in$ that operates at an `infinite speed'\footnote{Badiou, \textit{Being and Event}, 2006, p. 237.}. In modal logic, the symbol $\square$ conventionally designates the modality of \textit{necessity}---an operation that liberates logical determination from the instantaneous, totalizing constraints of set theory, making instead visible deduction as a process. In fact, this framework models not only the intuitionist logic of necessity, but also the dual modality of \textit{possibility}\footnote{The dual modality of possibility is formally derived as $\Diamond = \neg\square\neg$.}.

Yet, where Badiou presents this symbol of intervention $\square$ as an existentially discrete, atopological entity, the modal logic articulated by an elementary topos intersects with possible worlds semantics, formally acknowledging the structural multiplicity of different fidelities $\square$ as they originated in Saul Kripke's\footnote{Kripke, Saul (1981), \textit{Naming and Necessity}, Wiley--Blackwell.} work. Indeed, illustrating such a topological extension of Badiou's connection-operator, even in a pretopos world, a so-called Kripke-frame\footnote{A Kripke frame is defined as a pair $(W, R)$, where $W$ is a non-empty set of worlds and $R$ is a binary accessibility relation on $W$. For a given Kripke model $(W, R, \models)$ and a world $w \in W$, one writes $w \models \square A$ if and only if for all $u \in W$ satisfying $w R u$, it follows that $u \models A$; only under these conditions does the frame force the \textit{necessity} of the proposition $A$. The resulting modal system depends entirely on the geometric or categorical conditions imposed upon $R$.} contextualizes these suspended subjects of inquiry around a poset, directly extending it into a localic frame similar to Badiou's `transcendental indexing' $T$. 
Therefore, it is pivotal to ask how a topos envisions an entire multiplicity of such local operators, thus presenting a field connecting these very connections. 

Of course, in a general elementary topos $\Es$ there evidently exists no such simple external frame algebra. However, the subobject classifier $\Omega$ functions as an \textit{internal poset}, thus relationally generating the recursive rules of the so-called Kripke--Joyal semantics. In the sheafified architecture of a `weakly materialist' (Grothendieck-)topos, each  neighborhood $U$ within the site localizes a specific territory of possibility. 
Therefore, as I demonstrated above, Kripke--Joyal semantics effectively `territorializes' ostensibly disparate enquiries of  necessity, whereas within a local, atomistic topos of `democratic materialism', they are globally   reconciled into a single, unified intuitionist model of deduction graded over an external `transcendental indexing' $T$. It is only this externally bounded, `quasi-Platonic' theory of deduction that precludes torsion and satisfies Badiou's `fundamental law', endorsing his split theory of subjectivity as formulated in \textit{Being and Event}.

Inasmuch as Badiou's regulatory connection to an event thus acquires a structural mapping upon a predicate of necessity $\square$, one may conclude that a topos effectively bridges, in a more perpetual manner, \textit{mutually discrete fidelities} of the event. It constitutes a diagrammatic synthesis of distinct `enquiries', each executing its own chronological trajectory when topologically forced into a localic context of representation. Each enquiry, as Badiou presents it, is but a  \textit{localization} of the event---in fact, localization of \textit{`another event'} by which Badiou refers to the suspended process of inquiry itself, that is, a process \textit{enforced} by the split, `Platonic' outlook external to the topos in question. Symptomatically, the connection to Badiou's own presentation of the original event as $e$ remains entirely metaphoric in \BE. 

Therefore, even if a particular, `local' enquiry adopts a logico-\bk mathematical expression compatible with Badiou's formalism, the inquiry-event as a structural whole does not necessarily admit a global localizing frame in the explicit, elementwise manner that Badiou's `atomistic' perception of mathematical expression might suggest; it specifically does not have to collapse into a single localic topos of $T\text{-}\sets$. It is precisely the localic and reductive constraint in \textit{Logics of Worlds} that proves detrimental to mathematically understanding the expansive trajectory that the subject-philosophy of \textit{Being and Event} flirts with.

To evaluate this question of semantics and its materialization within the context of Cohen's procedure---which remains central to the subjectification traced in \textit{Being and Event}---one can interpret the Cohen topos as a \textit{diagrammatic} expression of the strategic position implicated by the procedure itself, rather than as the `generic', elemental outcome of that strategy (i.e. as $\female$). There is thus an inherent parallel between Badiou's presentation of Cohen's strategy and the urge to locate deductive fidelity into the formally disconnected regime, illustrating the shared `dialectical' motives driving both of Badiou's philosophical projects. 

Indeed, as I discussed, Cohen's construction itself fundamentally requires the strategic choice of a specific forcing poset $\Cond$ rather than the specific point $\female$---and it is precisely this structural selection that Badiou acknowledges as being `difficult', yet without engaging its philosophical core. It is precisely that structural choice that those in fidelity of general topos theory---those who have given up the `urge to Platonize'---need not atomistically capture; not even when endorsing the beauty of its expression as an internal effect (see Proposition~\ref{localeproposition}).

Therefore, the consequences that may be drawn from this evental rupture are far more complex than any specific Cohen topos $\Sh(\Cond, \neg\neg)$ can encapsulate. Topos theory arose precisely as a diagrammatic attempt to grasp the difficulties inherent to logical formalism, which by itself provides an inadequate frame for understanding its meta-ontological caesura. These difficulties, in fact, extend also beyond the strategic choice of $\cond$, let alone the `generic' point $\female$. Rather, they concern how subjectification itself in classical mathematics---unlike, say, in contemporary physics or geometry---occurs through the principle of divide and conquer, that is, the dialectical partialism still visible in Hegel's understanding of the `whole' through its \textit{parts}. That is precisely why the \textit{internal classifying object} $\Omega$ is called the `subobject classifier'; and that is precisely why an elementary topos, which was first discovered when Lawvere and Tierney adopted the atomistic topos of classical science as their object of study, is still \textit{internally atomistic}. 

Instead, what higher topos theory evidences is that in physics, as much as in mathematics, material reality does not generally admit the divisive, containing structure constitutive of the classical paradigm of science, ranging from Aristotle’s external logic to Hegel’s internalized conception of it. An elementary topos, by contrast, is still at risk of becoming `Platonically' self-contained through the language effectuated by its internal elements, in a fashion analogous to Badiou’s self-containing event, thereby supporting my characterization of the elementary topos as ultimately `singular’. Its singularity marks a condition that is not replicated in a higher topos, whose internal automorphisms make its truth `interactive’. A higher topos, too, admits a subobject classifier \(\Omega\), but this classifies only the so-called \((-1)\)-truncated arrows, whose fibres are either empty or contractible. The propositional, intuitionistic `dialect' associated with \(\Omega\) therefore fails to exhaust the homotopical structure of the topos as a whole---a fact that makes a higher topos categorically more \textit{difficult} than either Cohen or Badiou could envision.\footnote{See Veilahti, `Higher Theory and the Three Problems of Physics’, 2017.}

Obviously, mirroring Badiou's endeavour to escape the `Platonic' frame through Russell's paradox, there is a series of structural analogues illustrating the limits of internal atomism inherent to this `Hegelian' class of singular, internally dialectical topoi: similarly as the collection of different local models of sets (e.g. $\Sh(\cond, \neg\neg)$) do not form a localic topos, and the set of all sets is not a set, the collection of all elementary topoi does not form an elementary topos. Their collection is an expression of the \textit{plurality} rather than the singularity of truth-processes, as not all of them, under their internal formalism, can be synthesized as a `One-All' (as Badiouian commentators, such as Ray Brassier, might expect). But even if that `One-All' remains inconceivable as not merely a split, atomistic entity (Russell), but also as an elementary, internally `dialectical' one, category theory has already opened the Pandora's box of relational abstraction---a process that diminishes, step by step, the power of the divisive order of classification: first externally, and then internally, and God knows how category theorists will eventually reterritorialize categorical individuation. In short, the worlds envisioned by dialectical materialism do not themselves appear to form a single, `dialectical' entity---and this is the precise engine that, unlike what Badiou admits, drives the category-theoretic, `diagrammatic' abstraction as a process rather than an end point; something that demands the re-evaluation of all classical philosophy of being, from Parmenides to Heidegger. 

Likewise, the  event $e$, as it connects to local fidelities $\square_U$ or even the Kripke--Joyal semantics of a singular topos, is irreducible to any specific form of `totality'. This makes it a somewhat deceitful argument to say that an event, like the event of topos-theoretic expression, could only be locally intervened---through the `totality of its effects'. Instead, a topos \textit{per se} is an open, non-totalizing concept: it is a \textit{label for the specific and temporal processes of self-referentiality} that adopt multiple forms and implicate different relational layers, like topos theorists such as Lawvere and Tierney, and now also Jacob Lurie, have begun to study. At the same time, it deploys a process of mathematical innovation and abstraction, like the passage from Weil's incorporeal geometry to Grothendieck's relational insight, as a \textit{historical process internal} to mathematics. 
 Topos theory is thus a speculative, but at the same time `empirical' science precisely because it does not constitute a `statist' thesis or totality---it precisely is a process and `\textit{durée}'. 

This interpretation of the `difficulty' of forcing as an inherent, spontaneous feature of (scientific) subjectivity and its urge to always localize thus deviates from Badiou's `fundamental law' which externally enforces the primacy of a specific, `Platonic' model of forcing. If a topos is yet another event, it is the very \textit{modality of regulation} itself---and locality---which changes. Where Badiou's localic and atomistic concept of intervention merely recognizes the `statist' consequences of the event, there are now other means of inquiry which radically transform the regime of regulation itself: it no longer pertains to the split algebra of deductive logic, but expresses some deeper precursors capable of effectuating the very language carrying the weight of Badiou's dialectic event.

Consequently, as an altered regime of regulation---or rather their entire process---a topos serves as a structural analogue of Badiou's event (and \textit{durée}), engendering the apparent contradictions Cohen observed while eventually breaking the primacy of the atomistic notion of force. What a \textit{topos materially abstracts is the expression} of the problem of the `metaontological' or what, in a more Kantian vocabulary, accounts to the `transphenomenal real'.  It does not, purely speaking, abstract the meaning of what such categories designate independently of its own forms of expression; instead, it precisely articulates its own meaning as a relational process, thus counting as `another' event in its own right, and thus can only present itself to the `Platonic' minds through its formal, locally redundant effects.

As a result, Badiou's partial, atomistic cancellation of the distance between propositional enquiry and the effects of structural expression resonates with how he believes mathematics itself to fail to regulate, or express, a contingent and layered process: any sense of temporality could be only negatively endowed through the analytic incompleteness of mathematics. Thus, instead of grasping the plurality of the strategic reterritorializations enacted by genuinely diagrammatic modes of inquiry, Badiou's intervention remains faithful to 
his `narcissistic', non-interactional idea of the event: it can only be inquired by one subject at a time, and through \textit{only one order of necessity}---however temporal or suspended. Regulation itself 
remains predicatively split: a limitation preprogrammed to the very narrative describing the connection to fidelity $\square$. His idea of regulation remains narratively linear. As such, it is regulating necessity, rather than a necessary expression of mathematical regulation. 

Quite pertinently, this doctrine forecloses any option for his edifice to escape the centrifugal forces pressing his philosophy as it is trapped into circulating around the axiom of choice. It is `an axiom which', Badiou\footnote{Badiou, \textit{Being and Event}, 2006, p. 228.} assumes, `is intervention in-its-being', and thus, as I have demonstrated, structurally regulates the splitting of the classifying object as well as all epimorphisms, including the elementary arrow $\Omega \to 1$. 
This Badiouian approach to subjectivity and force, which are conceived purely in terms of splitting, then connects not so much to the event itself but to the precise doctrine that  precludes any temporality and eventality inherent to mathematical expression. It is \textit{dogmatic} specifically in the sense of structural expression itself being forced solely into the domain of logical constraints. There is no genuinely `spontaneist' fidelity in mathematics, according to Badiou, which ultimately renders his own engagement with Kripke-algebras and localic topos theory ironic. 

The analogue I have proposed between the event, which is non-local, and topos-theory, which reflects a specific encounter between global expression and local representation, is undoubtedly `unassignable' to the Badiouian structure of enquiry, which is formally imprisoned by Democritus' problem of parts. Anything novel about a topos is read strictly from within this framework, as it forms `a particularly nonsensical part' of mathematics\footnote{Ibid., p. 237.}. Therefore, to Badiou, the only sensible form of regulation is the logical, incorporeal algebra of necessity, blinding him from the contingency of relational expression.

This dogmatism, which denies the role of inscription (also as a matter relevant to Badiou's own philosophy), is visible directly in how these inquiries are constrained in Badiou's own writing, as he never questions the isolated enquiries in \textit{Being and Event} nor the `democratic' or `atomistic' condition characterizing topoi in the case of \textit{Logics of Worlds}. Badiou\footnote{Ibid., p. 233.} thus maintains that `every fidelity is particular', but consequently does not take the `difficult' step of acknowledging that its expression, and not only what it designates, puts at risk `the universal \textit{form} of the procedures which constitute it'. Thus, to Badiou, events themselves become discrete and isolated, like the `isolate' in his own `calculated phenomenology': an event might exist independently of the `inexistent' or \textit{nihil privatum} of `Platonic' idealism, but nonetheless its voice remains muted by the strict, `democratic' topology of language. 
The event cannot then but incur an infinitude of isolated subjects and histories, each articulating itself as the global truth; this particularly prevents him from understanding continuities and interactions between several events, let alone multiple perspectives on their temporal expression. 

Through a topos, by contrast, that which is `unassignable' in Badiou's philosophy, and thus carries a metaphoric name, manifests itself precisely \textit{not} as `One-All', or as merely a `caesura'. Rather, it emerges as time itself and an abstract process that alters the very concept of expression---a process that is a body, even if not in a singular, isolated sense of idealism.

\section{Topos: an Event, Change, and a Body}

The relational `essence' of topos theory is structurally distinct from the logical essence of set theory, both regarding their forms of expression as well as their internal designation of truth. Insofar as topos theory constitutes a ` ``universe of discourse'' within which we can carry constructions'\footnote{Johnstone, \textit{Topos Theory}, 1977, p. 47.}, it topologizes not solely the constructions themselves, but also how it perceives its own discourse, that is, its internal concept of such very constructions. It localizes the very \textit{discourse of being}---perhaps more so than the image of being this incorporeal discourse prescribes. 

If Michel Foucault\footnote{Foucault, Michel (2002), \textit{The Archaeology of Knowledge.} Trans. A. M. Sheridan Smith. London and New York: Routledge [Originally published in 1969.]} refers to this discursive register of knowledge as an `archive', topos theory demonstrates that it is the archive itself---not that which is archived---which emerges \textit{topo}-logically, alongside the processes pertinent to the formation of space. This is precisely because any `totality of its effects'\footnote{Sartre, Jean-Paul (2003), \textit{Being and Nothingness. An essay on phenomenological ontology}, London and New York: Routledge. p. 1.}, which Foucault conceives as an ultimately `archaeological phenomenon', cannot itself be archived within a classical hierarchical schema, any better than the 1-categorical schema of elementary topos theory. 

Therefore, where classical mathematics conceives torsion as something obscuring the hierarchical, `Platonic' doctrine of the archive, topos theory poses a different, more relational approach to archives: they are not so much the starting point of mathematical inquiry, but precisely an effective outcome of the more or less \textit{temporal} processes of expression. 

A more conservative philosopher of mathematics might now align with Badiou's tenet that, despite the `new' modes of expression topos theory performs, it is nonetheless `the same mathematics' that dominates the field. Such an assertion is clearly not far-fetched with respect to how a mathematician still makes points about categories---one still aims at propositions, even about the apparently altered expression of the scientific archive. By imposing localizing functors, however infinite, one still \textit{forces} logical claims, even about the relational, internally displaced picture of being. 

Badiou might thus declare his victory, arguing that topos theory has at best modified how one charts the truth as a process; perhaps the tools through which we map being are ostensibly different, but that `being' itself which is the ultimate landscape covered by these maps remains intact. 

Yet, at least the procedural shift is undeniable: these maps represent different, locally effective `snapshots' of truth without accounting for the situation as an integrated whole. This framework suggests an alternative account of the materialization of knowledge which fundamentally falsifies any epistemic theory grounded upon localic, hierarchical poset-structures; it falsifies the \textit{necessary} organization of knowledge according to binary points just as it challenges the external transcendental structures pertinent to Badiou's `calculated phenomenology' of `democratic materialism'. In short, it falsifies a `constitutive' order of scientific materialism, despite mathematicians' attempts to build a new propositional hierarchy and recapture the new `forms of sensibility' in an altered, topological way. Knowledge and science, as they emerge as part of their own representation, are no longer asymptotic in the way C. S. Peirce's philosophy of science naïvely suggests, despite all the effort to reclaim a new `archimedean' insight enacted at the procedural level of proofs. 

The main argument of this paper is that, despite these continuing efforts to convey `atomistic' points about this increasingly geometric `archive' of knowledge, Badiou's philosophy as such does not preclude the claim that \textit{topos-theoretic expression itself counts as an event}. Perhaps it is not an isolated, merely self-referential event-multiple $e \in e$ formalizing how Badiou himself pursues to capture its local consequences in set-theoretic domain. Likewise, I do not argue that any particular topos $\Es$, as an ontological, representable entity, could account for the abstract event holding its sensibility and formalism together---a concrete instance of a topos is not identical to its scientific `body'. 
Instead, a topos, nonetheless, is something that genuinely \textit{happens} to this body. What truly changes through a topos is how it alters the question of topology and the very form of the mathematical `letting-be-of-appearing': ultimately the `specific manner[s] of holding-together'\footnote{Badiou, \textit{Being and Event}, 2006, p. 127.}.

Let me put it in this way: if an event in Badiou's philosophy can be only accessed through its \textit{local consequences}, this, indeed, forms the crux of the `gap' that separates the event nominated `from above' and the assessment of its local consequences `from below'. 
But now, through topos theory, one acquires a topological and in some sense increasingly temporal chart of this `gap' situated amidst scientific expression and the propositions and points through which it is procedurally archived. In fact, topos theory can entirely bridge this gap by declaring how the event, `from above', actually happens to the very concept of what it is to enquire its being `from below'. 

Hence, my argument does not contradict Badiou's core hypothesis of the dialectical conflict between the event, perceived as an effective and nominal totality, and its local gestures visible in the archives. Quite the contrary, one truly in fidelity of the event must acknowledge that it is mandatory also for the meaning of the archive, or `local configuration', to change---that, indeed, is the essence of topos theory! It is the making visible, or assigning otherwise that what formerly remained `inexistent' or `unassignable'---the variety of options for the subject of `forcing', but also a multiplicity that concerns subjectivation (and topologization) as concepts, that is, topos theory alters the very \textit{act} of forcing. This `event', or change in subjective expression, is thus structurally maximal: it is through this `retroaction'\footnote{Ibid., p. 452.} enacted by a topos, as it questions its own composition, that `the guarantee of a certain unity in being, beyond its surge into appearing', acquires a genuine `corporeal synthesis'.

Badiou's exposition, unable to capture this structural mutation of topology and `subjective formalism' itself, is led into a mere paradox, asserting that the scientific effects of `[t]he immanent agent of the consequences (of a truth), or the possible agents of their denial, or that which renders their occultation possible', are all bounded by the old, logical archive\footnote{Ibid., p. 453.}.

Instead, my claim is that the `local configuration' itself is not a `statist' or atomistic world, like it appears in Badiou's symptomatic claim towards an ultimate, fixed expression of event-philosophy. While Badiou's concept of local subjectivity remains cornered into a logical narrative---his only metaphor of formal inscription, presented as `the general legislation of appearing' constitutive to any mathematically conceivable `phenomenology'---outside this atomistic peninsula of `democratic materialism', change itself needs not be treated only through an isolated object. Change needs not be treated merely as a question of how an `object can ``become'' a site' in some preprogrammed, even if `paradoxic' way---a type of change that, when atomistically observed, inevitably `vanish[es] without delay', denying the true `amplitude of its consequences'\footnote{Ibid., p. 452.}.

Instead, for an event to truly appear as evocatively as Badiou claims, the `democratic' concept of science needs a critical re-evaluation, as does the question of how mathematics itself, in a formally tangible way, traces structures analogous to those of dialectical materialism, and even beyond. Badiou offers a world split between two, mutually exclusive alternatives: either `democratic' or `dialectical'. This precisely explains why \textit{scientific change}\bk ---the fundamental paradigm shift from set theory to category theory---proves so intractable for Badiou to articulate. 

From another angle, this profound shift is often also missed by `constructivist' approaches to science studies, despite their otherwise open definition of social categories. For instance, Andrew Pickering's\footnote{Pickering, Andrew (1984), \textit{Constructing Quarks: A Sociological History of Particle Physics}. Chicago: University of Chicago Press.} and Peter Galison's\footnote{Galison, Peter (1987), \textit{How Experiments End}. Chicago: University of Chicago Press.} elaborations of physics come to mind, as materialism is still presented as something external to mathematical reasoning, instead of viewing change as a concept entangled with the very notion of formal expression. The `material practices' referred to by conventional science studies thus fundamentally preclude the domain of mathematics, endorsing a split separation between theory and empiria, as opposed to the more poststructuralist idea of speculative empiricism, where mathematical structuralism\footnote{See Awodey, Steve (1996), `Structure in Mathematics and Logic: A Categorical Perspective', \textit{Philosophia Mathematica}, 4(3), pp. 209--237; Shapiro, Stewart (1996), `Mathematical Structuralism', \textit{Philosophia Mathematica}, 4(2), pp. 81--82; and Shapiro, Stewart (2005), `Categories, Structures, and the Frege--Hilbert Controversy', \textit{Philosophia Mathematica}, 13(1), pp. 61--77.} itself emerges as a product of experimentation.

Unlike previous accounts exploring the philosophy of category theory by focusing strictly on the structuralist organization of mathematics, I thus interrogate whether topos theory can provide an alternative insight into the \textit{post}structuralist orientations of the `French moment' of philosophy, to use Badiou's own phrase, including profound re-readings of dialectics and materialism. In this sense, the picture science studies offers is inescapably secular. Poststructuralism, instead, is a tradition that responds to dilemmas inherent to such classical structural narratives---those specific points, or local configurations of structure, which Badiou stubbornly frames as logical rather than geometric. 

Badiou's `dialectical materialism' aligns with this task as it seeks an escape 
route away from the localic peninsula of classical scientific formalism. Yet, its 
success remains incomplete, for he believes that science itself can offer at best 
a novel `fact', if not merely a `modification'---a genuine change, instead, 
requires something `singular' according to Badiou\footnote{Badiou, \LW, 2009,   
pp. 47, 58--61.}; something science itself could only recognize as an `occultation 
of truth'. My claim, here, is that topos theory is a double-edged sword. It both 
expresses the limitations inherent to the atomistic representation of truth 
pertinent to localic theory, thus offering a formal articulation of that precise 
path \LW\ was obsessed with finding. Yet, by charting this philosophical passage, 
topos theory also declares this `singular' expression of change all too 
profane---an elementary topos, after all, remains internally atomistic, resembling 
the dialectical concept of an `artistic truth [where] the world exhibits a 
singular form of tension between the intensity of the sensible and the 
tranquillity of form'. Change, then, is articulated solely in an elemental way, 
as a matter for an individual body. It is, by all means, an important stepping 
stone in bridging the philosophical partialism implicated by the atomistic peninsula with the genuinely relational and 
interactional understanding of temporality. But singular change, on its own, 
fails to articulate the movement, or `\textit{durée}', between atomistic logic and 
the relationally \textit{higher} structures of science. 
Elementary topos theory, as this paper has framed it, is thus precisely a 
bridge---not an end in itself.

\vspace{.5cm}
\noindent \textbf{Topos-theoretic Typology of Change}
\vspace{.2cm}

To formalize this bridge in terms of elementary topos theory, Badiou's concept of change---when applied to mathematical expression itself---indeed acquires a precise, topos-theoretic interpretation: it yields a genuinely new model of the production of logical expression; not just a fact but a rupture that demands every claim about the world to be structurally re-evaluated. 
Cohen's procedure, beginning from the original, set-theoretic situation $\sets$, and evolving through the factual emergence of new points all the way to the diagrammatic multiplicity of the strategic, logically contradicting effects endorsed by each individual $\cond$, articulates a precise narrative of how the existence of a topos becomes `maximal' in relation to the old, localic  body of mathematical expression. 
 The crux of this passage is precisely the double articulation of a topos as both a category of understanding (topos) and as that dialectical 
language of sensibility (the Mitchell--Bénabou language)---a bridge that permits 
Badiou's four-fold conception of change to be situated within the topos-theoretic 
`story of the event'.

\subsubsection*{I. Modifications}

There exist different models of set theory ($\Ss$) constructed during Cohen's 
reformulated forcing demonstration examined above. A Cohen topos is not strictly two-valued, like the `Platonic' situation, but it is still only a `modification' because it can be forced into a two-valued situation by the ultra-power construction, and thus nothing has truly changed: it contains no new, 
previously inexistent structures; the void remains as 
empty as it ever was, and its internal logic endorses the law of the 
excluded middle.

\subsubsection*{II. Facts}
Instead, a non-Boolean locale, like the general form of Badiou's world or $T$-sets, or sheaf theoretically $\Sh(T)$,  can produce a genuine `fact' if it does not give 
rise to the axiom of choice (AC). As a `quasi-split' situation, nonetheless, it is not a true singularity, for the intuitionists the axiom of choice itself was never certain to begin with. Any `singularity' or torsion 
that could render its internal semantic truly consequential to the 
`fundamental law of the subject' remains absent. In this reactive reading of topos theory compatible with Badiou's formalism, its expression as such does not engender new points: the ambiguity Cohen recognized is an inherent aspect of set-theory itself, even if proven consistent with the categorical formalization---structures understood in the `second sense', to borrow a Kantian phrase, like this paper evidences. Its contributions to those still `wishing to Platonize' remain a matter of belief.

\subsubsection*{III. Weak Singularities}
In contrast, a Grothendieck topos, which my previous paper labeled as `weakly' or `dialectically materialist'\footnote{Veilahti, \textit{Alain Badiou's Mistake}, 2013.}, falsifies Badiou's 
strong, `atomistic postulate' of classical materialism. In other words, that world offers no `real synthesis' but its more `corporeal synthesis', as enacted by sheaf theory, may exhibit torsion that genuinely contradicts the `fundamental law', as such a site is not externally graded by a neutralizing `transcendental indexing' like Badiou's Heyting structure $T$. In fact, within such a topos, there may exist multiple 
distinct `corporeal' or material syntheses $\Es \to \sets$, which are not structurally unique precisely when it is not a localic world (i.e. one where $\Omega$ is globally generated, or SG), that is, precisely when it falsifies Badiou's logically enforced definition of the `coherence of bodies'\footnote{In Badiou's framework, `\textit{a body $C_{\epsilon}$ relative to the proper inexistent $\epsilon$ of a site $(A, \text{Id}_A)$, is the set of elements of $A$ which incorporate themselves to the evental present}' of the site $A$. These points at which bodies adhere are conditioned by the principle of the `\textit{coherence of bodies}', which draws from the existentially governed, idealist transcendence $T$. This structure ensures that the bodies `admit[s] a real synthesis, that is of an envelope for the order-relation $\leq$, which operates not on phenomena but on the multiple composition'---for the element $S_c \in A$ the body $C_{\epsilon} = \downarrow\!(S_c)$ is reflected as the `corporeal synthesis'. Badiou, \textit{LW}, 2009, pp. 485--486.}. It thus contradicts the \textit{Tractatus} and the factual 
existence `of a single world'\footnote{Ibid., p. 358.}, but does not yet render its existence strongly `singular' 
or immaterial---its expression becoming \textit{material} in terms of classical formalism of $\sets$ is still an option. Namely, as Badiou\footnote{Ibid., p. 487.} declares, `unless the inexistent has become over-existent, the elements of a post-evental body do not maximally exist in the site. They have to become'.

\subsubsection*{IV. Strong Singularities}
\noindent Only the general class of elementary topoi credits such sites of `becoming' that 
may be regarded as `strong singularities'---sites truly occulting the `Platonic', 
idealist expression of set-theoretic ontology. Such a topos does not `materialize' over the atomistic topos of sets even through its sheafification. Points, if they exist, thus possess a 
different, more geometric meaning: they can express non-trivial \textit{automorphisms}, which Badiou's localic formalism forcefully neutralizes. By recognizing that the `local' (or topological) is not the same as localic or atomistic, an elementary, or `singular' topos opens up a surface that can be read precisely as a situation from which there emerge `artistic subjects'\footnote{Badiou, \BE, 2006, p. 73.}. These are subjects of an `artistic truth [where] the world exhibits a singular form of tension between the intensity of the sensible and the tranquillity of form'.

This makes an elementary, non-Grothendieck topos `strong' or maximal with respect to its consequences over the localic site. In such a site, the classical `forms of sensibility' are fully internal effects of the world---it is precisely a site of \textit{becoming}, like the existence of the Mitchell--Bénabou-language and Kripke--Joyal semantics demonstrate: set-theoretic expression adopts purely `categorical understanding'. A topos is literally the `bearer of the subjective appearance of truth'\footnote{Badiou, \textit{Logics of Worlds}, 2009, p. 451.} from which the classical, intuitionist language of sets is `borne': `\textit{[a] multiple-being which bears this subjective formalism and thereby makes it appear in a world}', thus receiving `\textit{the name of ``body''}'\footnote{Ibid., p. 453.}.

\vspace{.5cm}

This bridge between topos theory and the typology of change allows a more evocative understanding of Badiou's concept of a body. 
Yet, again, two observations can be drawn from this association. First, a poststructuralist account must conceptualize its semantic endowments as an \textit{effective} result of the relational rearticulation of forcing operations. Language is not only an analytically unnecessary, or `synthetic' category, but forcing itself, as a point of subjectivation, becomes a \textit{differential concept}, expressing a variety of layers of relational structure. After all, even Badiou must admit this: `post-evental bod[ies]' cannot exist in the site, but `have to become'.

But how can the concept of a body truly intersect with mathematics? 
To Badiou himself, of course, mathematics proper remains strictly exiled from this `poetic' register of the event. Yet, his argument is inextricably meta-ontological, whereas my question concerns mathematics as \textit{expression}---precisely not its `being' or what it `designates'. Just as Badiou argues that Paul Val\'ery's famous poem, `The Graveyard by the Sea' (\textit{Le Cimeti\`ere marin}), is `a story of an event', mathematics, insofar as it adopts new modes of expression, is a similar story. The internal classifier object $\Omega$ marks that `undivided unity' of its expression that the classical, split language of set theory failed to envision. As an abstract concept, $\Omega$ is a category-theoretic `manner of holding-together'\footnote{Badiou, \textit{Being and Event}, 2006, p. 127.} the local faces of `subjective formalism'; it is like `the sea' in Val\'ery's poem that constitutes `the surface of the world', whereas `the sun' is then like the terminal object---`that facet of the world' acting as the `fixed point' of the sea; a confirmation or `letting-be of appearing' of its `undivided unity'\footnote{Badiou, \LW, 2009, p. 455.}. Moreover, in a `weakly material' topos, which in the aforementioned typology is at the same time `weakly singular', this `sun', then presented as a sheaf, provides `the function of appearing of the sea', precisely linking it to a  `transcendental functor' that is Badiou's own name for a sheaf.

Badiou's concept of a body is also irrevocably existential. And an elementary, `singular' topos inevitably expresses a sort of `consciousness' about its own, internal formalism---its own `thought of $\mathcal{N}\text{ature}$'. Insofar as Badiou's disbelief in the radical effectivity of topos theory is maintained, this consciousness certainly remains one which `is but the vain shadow of Noon and that the dead, ensconced in the earth, are the symbolic future of this ``change''';\footnote{Ibid., p. 456.} a mutation that becomes cancelled `at the hands of Parmenides and Zeno, at the hands of those who refuted movement'\footnote{Ibid., p. 457.}. And, ironically, the localic, `calculated phenomenology' of \LW\ precisely carries the weight of this refutation. It perfectly fits a narrative where, `[i]n this entire affirmative peroration, the [topos-theoretic] ``body'' noticeably gathers together disparate elements, all of them drawn from the storm rising over the sea: ``the wind's reviving'', the ``salt-breathing potency'', [\dots]'\footnote{Ibid., p. 466.} These are all effects of the `dialectical' attempts to regain control by enforcing a sheafified, `Platonic' form of materiality over the topos.

Therefore, providing a structural narrative of mathematics itself as a process, topos theory articulates \textit{novelty} and change in a way expanding in two different directions---both appealing to those admiring the old, `dialectical' reaction to the logical peninsula of science, but also to those favoring more open, poststructuralist orientations. 
This particularly relates to the second observation stemming from the link between Badiou's typology of change and topos theory: the fact that an elementary topos, while acquiring a categorical body that endorses Kant's transphenomenal processes in the `second sense', still operates in the domain of $1$-categories, appealing to the singular orientation of event-philosophy. 

It is, indeed, a bridge, or the first layer of the `differential', poststructuralist process amidst the logical and categorical frameworks, but still only an internally `dialectical' expression consistent with Badiou's theory of change as \textit{a} body. The central weakness of the Badiouian project stems from his attempt to theorize the body, and thus also that what is set to `become', as an isolated `object'. A body, to him, is like a single site or topos that `vanishes without delay', despite `its consequences sometimes characteriz[ing] it as an event'\footnote{Ibid., p. 452.}. It is particularly not an \textit{interactive} and \textit{temporal} process between different theories of a topos, whereby the very structures through which change is embodied are permanently at question.

What might then be a kind of philosophical, but also political thinking that discards a static axiomatic, instead conceiving situations and events themselves differentially and processually, following the motives behind Gro\-thendieck's and Lurie's insight? If the answer lies somewhere in the vicinity of poststructuralism, such thinking may, at the very least, recognize the prevalence of the virtual itself. Even if Badiou sides with logical reactionaries, arguing that the consequences of category-theoretic expression are nil, and that nothing new has actually been established, the question is precisely about how topos theory reframes the relationship between novelty and actuality. Even if nothing is actually new, novelty and difference can be viewed as `virtual' or geometric dimensions of being itself---something not intrinsically coordinated by the conditions of actuality that any specific axiomatic demands. 

Topos theory---whether or not it actualizes that difference---is an expression of this virtual plane; it illustrates how different tools of actualization, or axiomatics, contribute to a whole variety of understandings of points, separate from the classical, propositional structure of science. Undoubtedly, Badiou's own politics of mathematics remains belated, refusing to endorse such developments. Yet it is precisely when one moves to a `higher' topos, beyond the singular framework of elementary topoi, that the gap between Grothendieck's summit and Badiou's logical peninsula becomes bridged, opening up an entire plane of expressions in which difference can individually actualize.

\section*{References}



\setlength{\parindent}{0pt}

\setlength{\parskip}{1.5ex}
\small{

Anel, Mathieu \& Gabriel Catren (eds.) (2022), \textit{New Spaces in Mathematics: Formal and Conceptual Reflections}, Cambridge University Press. 

Awodey, S. (1996),  `Structure in Mathematics and Logic: A Categorical Perspective'. \textit{Philosophia Mathematica} 4 (3). pp. 209--237. doi: 10.1093/ philmat/4.3.209.

Badiou, Alain (2006), \textit{Being and Event.} Transl. Oliver Feltham. London, New York: Continuum. [Originally published in 1988.]

Badiou, Alain (2009), \textit{Logics of Worlds. Being and Event, 2.} Transl. Alberto Toscano. London and New York: Continuum. [Originally published in 2006.]

Badiou, Alain (2022), \textit{The Immanence of Truths: Being and Event III}, translated by A. Toscano and K. Reinhard, Bloomsbury Academic.

Bowker, Geoffrey C. \& Susan Leigh Star (1999), \textit{Sorting Things Out: Classification and its
Consequences.} Cambridge: MIT Press.

Cohen, Paul J. (1963),  `The Independence of the Continuum Hypothesis', \textit{Proc. Natl. Acad. Sci. USA} 50(6). pp. 1143--1148. 

Cohen, Paul J. (1964),  `The Independence of the Continuum Hypothesis II', \textit{Proc. Natl. Acad. Sci. USA} 51(1). pp. 105--110.

Deleuze, Gilles (1994), \textit{Difference and Repetition}. Transl. Paul Patton, New York: Columbia University Press. [Originally published in 1968.]

Deleuze, Gilles \& Félix Guattari (1988), \textit{Thousand Plateaus---Capitalism and Schi\-zophrenia}. Transl. Brian Massumi, Minneapolis, University of Minnesota Press. [Originally published in 1980.]

Foucault, Michel (2002), \textit{The Archaeology of Knowledge.} Trans. A. M. Sheridan Smith. London and New York: Routledge. [Originally published in 1969.]

Fourman, Michael P. (1980),  `Sheaf models for set theory', \textit{Journal of Pure and Applied Algebra} 19. 91--101.

Galison, Peter (1987), \textit{How Experiments End}. Chicago: University of Chicago Press.

Johnstone, Peter T. (1977), \textit{Topos Theory}. London: Academic Press. 

Johnstone, Peter T. (2002), \textit{Sketches of an Elephant. A Topos Theory Compendium.} Volume 1. Oxford: Clarendon Press. 

Kripke, Saul (1981), \textit{Naming and Necessity}, Wiley--Blackwell.

Krömer, Ralf (2007), \textit{Tool and Object: A History and Philosophy of Category Theory}. Science Networks. Historical Studies 32. Berlin: Birkhäuser.

Mac Lane, Saunders \& Ieke Moerdijk (1992), \textit{Sheaves in Geometry and Logic. A First Introduction to Topos Theory}. New York: Springer-Verlag.

Marquis, Jean-Pierre (2009), \textit{From a Geometrical Point of View. A Study of the History and Philosophy of Category Theory.} New York: Springer.

Nietzsche, Friedrich (2002), \textit{Beyond Good and Evil}. Transl. Judith Norman. Cambridge: Cambridge University Press.

Pickering, Andrew (1984), \textit{Constructing Quarks: A Sociological History of Particle Physics}. Chicago: University of Chicago Press. 

Sartre, Jean-Paul (2003), \textit{Being and Nothingness. An essay on phenomenological ontology}, London and New York: Routledge.

Shapiro, S. (1996),  `Mathematical structuralism'. \textit{Philosophia Mathematica} 4(2), 81--82.

Shapiro, S. (2005)  `Categories, structures, and the Frege-Hilbert controversy: The status of meta-mathematics'. \textit{Philosophia Mathematica} 13(1), 61--77.

Veilahti, Antti (2013), `Alain Badiou's Mistake: Two Postulates of Dialectical Materialism'. \textit{Preprint}, arXiv. Version 4 revised in 2026. \url{https://arxiv.org/pdf/1301.1203}. 

Veilahti, Antti (2017), `Higher Theory and the Three Problems of Physics', \textit{Preprint}, arXiv. \url{https://arxiv.org/pdf/1712.09454}.

Wittgenstein, Ludwig (1981), \textit{Tractatus Logico Philosophicus}. Transl. C. K. Ogden. Oxon and New York: Routledge. [Orig. published 1922.]

Zalamea, Fernando, \textit{Synthetic Philosophy of Contemporary Mathematics}, translated by Z. L. Fraser, Urbanomic, 2014.

}

\end{document}